\documentclass[12pt,a4paper]{amsart}
\usepackage{amsmath,amsthm,upref}
\usepackage[dvips]{graphicx}
\usepackage{psfrag}
\usepackage[all]{xy}
\usepackage{a4wide}
\usepackage{epsfig,subfigure}
\usepackage[hypertex]{hyperref}
\usepackage{color}
\usepackage{a4wide}
\numberwithin{equation}{section}

\newcommand\Sl{\textrm{SL}_2}

\newcommand{\R}{\mathbb{R}}
\renewcommand{\P}{\mathbb{P}}
\newcommand{\C}{\mathbb{C}}
\newcommand{\Q}{\mathbb{Q}}
\newcommand{\Z}{\mathbb{Z}}
\newcommand{\N}{\mathbb{N}}
\newcommand{\CCC}{\mathcal{C}}

\newcommand{\QQQ}{\mathcal{Q}}

\newcommand{\XXX}{\mathcal{X}}
\newcommand{\HHH}{\mathcal{H}}

\newcommand{\LLL}{\mathcal{L}}
\newcommand{\comb}{\mathcal G}
\newcommand{\marking}{m}
\newcommand{\SL}{{\rm SL}}
\newcommand{\GL}{{\rm GL}}

\newcommand{\ol}{\overline}

\newcommand{\frach}{\mathfrak{h}}
\newcommand{\fracg}{\mathfrak{g}}
\newcommand{\fraca}{\mathfrak{a}}
\newcommand{\fracu}{\mathfrak{u}}
\newcommand{\fracn}{\mathfrak{n}}

\newtheorem{Theorem}{Theorem}[section]
\newtheorem{Corollary}[Theorem]{Corollary}
\newtheorem{Lemma}[Theorem]{Lemma}
\newtheorem{Claim}[Theorem]{Claim}

\theoremstyle{remark}
\newtheorem{Remark}{Remark}[section]

\theoremstyle{definition}
\newtheorem*{Acknowledgments}{Acknowledgments}

\newtheorem*{Convention}{Convention}

\begin{document}
\title[The Arnoux-Yoccoz Teichm\"uller disc]
{The Arnoux-Yoccoz Teichm\"uller disc}

\author{Pascal Hubert, Erwan Lanneau, Martin M\"oller}

\address{
Laboratoire d'Analyse, Topologie et Probabilit\'es (LATP) \newline
Case cour A Facult\'e de Saint J\'er\^ome \newline
Avenue Escadrille Normandie-Niemen \newline
13397, Marseille cedex 20, France
}

\email{hubert@cmi.univ-mrs.fr}

\address{
Centre de Physique Th\'eorique (CPT), UMR CNRS 6207 \newline
Universit\'e du Sud Toulon-Var and \newline
F\'ed\'eration de Recherches des Unit\'es de 
Math\'ematiques de Marseille \newline
Luminy, Case 907, F-13288 Marseille Cedex 9, France
}

\email{lanneau@cpt.univ-mrs.fr}

\address{
Max-Planck-Institut f\"ur Mathematik\newline
Postfach 7280\newline
53072 Bonn, Germany
}

\email{moeller@mpim-bonn.mpg.de}

\date{\today}

\begin{abstract}
We prove that the Teichm\"uller disc stabilized by the Arnoux-Yoccoz 
pseudo-Anosov diffeomorphism contains at least two closed Teichm\"uller geodesics. This
proves that the corresponding flat surface does not have a cyclic Veech group.

In addition, we prove that this Teichm\"uller disc is dense inside the
hyperelliptic locus of the connected component $\mathcal H^\textrm{odd}(2,2)$. 
The proof uses Ratner's theorems. 

Rephrasing our results in terms of quadratic differentials, we 
show that there exists a holomorphic quadratic differential, on a genus 
$2$ surface, with the two following properties. \\
(1) The Teichm\"uller disc is dense inside the moduli space of holomorphic quadratic differentials 
(which are not the global square of any Abelian differentials). \\
(2) The stabilizer of the $\textrm{PSL}_2(\R)$-action contains two non-commuting 
pseudo-Anosov diffeomorphisms.

\end{abstract}

\maketitle
\setcounter{tocdepth}{1}
\tableofcontents

\section{Introduction}

After the work of McMullen~\cite{Mc1,Mc2, Mc4,Mc5,Mc6} and Calta~\cite{C} 
who classified Veech surfaces and completely periodic translation surfaces 
arising from Abelian differentials in genus 2, it is a big challenge to 
try to understand what happens for translation surfaces that admit
pseudo-Anosov diffeomorphism in higher genus. We recall that the trace field of a pseudo-Anosov 
diffeomorphism is $\Q[\lambda + \lambda^{-1}]$ where $\lambda$ is the expansion factor. This is an 
invariant of the underlying translation surface: it is the holonomy field. The trace field 
of such a surface is 
a number field of degree bounded above by the genus of the surface
(see Appendix of Kenyon-Smillie \cite{KS}). At the time
Thurston defined pseudo-Anosov diffeomorphism, no examples were
known when the holonomy field is an extension of $\Q$ of degree more than 2.

At the beginning of the 1980's, Arnoux and Yoccoz~\cite{Arnoux:Yoccoz} 
discovered a family
$\Phi_n$, $n\geq 3$ of pseudo-Anosov diffeomorphisms with
expansion factor $\lambda_n=\lambda(\Phi_n)$ the Pisot root of the
irreducible polynomial $X^n-X^{n-1}-\dots-X-1$. The pseudo-Anosov diffeomorphism
$\Phi_n$ acts linearly on a genus $n$ translation surface $(X_n, \omega_n)$ 
(the Abelian differential $\omega_n$ having two zeroes of order $n-1$). 
The stable and unstable foliations of these pseudo-Anosov diffeomorphisms exhibit 
 interesting properties studied in~\cite{Ar1}. See also \cite{Ar2} 
where a strange phenomenon is discussed in the case of genus $3$.
 
In this paper, we study properties of the Teichm\"uller disc 
stabilized by the pseudo-Anosov diffeomorphism $\Phi = \Phi_3$. The 
Abelian differential $\omega = \omega_3$ has two zeroes of order
$2$. The stratum of Abelian differentials (see Section~\ref{strata} 
for definitions) with two zeroes of order $2$ is called $\HHH(2,2)$. 
It has two connected components (see Kontsevich-Zorich~\cite{KZ}). The
translation surface $(X,\omega)$ sits in 
$\HHH^{odd}(2,2)$ the non-hyperelliptic component with odd parity of
the spin structure. Nevertheless, as $\Phi$  
was first defined on the sphere and then lifted to $(X,\omega)$, the 
Riemann surface $X$  is a hyperelliptic surface.
Thus, $(X,\omega)$ belongs to the hyperelliptic locus in
$\HHH^{odd}(2,2)$. We denote this locus by $\LLL$. It has 
complex codimension $1$ in $\HHH^{odd}(2,2)$. 

Hubert and Lanneau~\cite{HL} showed that this Teichm\"uller disc was not
stabilized by any parabolic element. Thus 
the Arnoux-Yoccoz surface has been considered as a good candidate
for a Teichm\"uller disc stabilized only by a cyclic group
generated by a pseudo-Anosov diffeomorphism. We show that this
does not hold.

\begin{Theorem}\label{pseudo-Anosov}
The surface $(X,\omega)$ has at least two transverse hyperbolic directions. 
More precisely the Teichm\"uller disc of $(X,\omega)$ is stabilized by 
$\Phi$ and by a pseudo-Anosov diffeomorphism $\tilde \Phi$ not 
commuting with $\Phi$. The expansion factor 
$\tilde \lambda=\lambda(\tilde \Phi)$ has degree $6$ over $\mathbb Q$.
\end{Theorem}

We recall that the stabilizer of a translation surface under the 
natural $\Sl(\R)$-action is a Fuchsian group called the Veech group. 
The Veech group is the group of the differentials of the diffeomorphisms 
that preserve the affine structure induced by the translation structure.
Combining Theorem~\ref{pseudo-Anosov} with the aforementioned result of 
Hubert-Lanneau (\cite{HL}) we get the following corollary.

\begin{Corollary}
Non-elementary Veech groups without parabolic elements do exist.
\end{Corollary}

\begin{Theorem} 
\label{orbit-closure}
The Teichm\"uller disc stabilized by the Arnoux-Yoccoz
pseudo-Anosov diffeomorphism $\Phi = \Phi_3$ is dense inside the hyperelliptic locus
$\LLL$.
\end{Theorem}
 
Theorem~\ref{pseudo-Anosov} and Theorem~\ref{orbit-closure} show that 
in genus $3$ the situation is very different from what we know in 
genus $2$. McMullen 
proved that, as soon as a genus 2 translation surface is stabilized by 
a pseudo-Anosov diffeomorphism with orientable stable and unstable
foliation, the image of its Teichm\"uller disc in the moduli space
of curves is contained in a Hilbert modular 
surface. In particular, it is not dense in its stratum. 
Here, the Veech group is non-elementary, nevertheless, 
the closure of the Teichm\"uller disc is as big as it can be. 

It is also interesting to rephrase Theorems~\ref{pseudo-Anosov}
and~\ref{orbit-closure} using quadratic differentials. In particular
we show that the orientability assumption in genus $2$ is necessary.

\begin{Corollary} 
\label{quadratic-diff1}
Let $(\P^1, q)$ be the quotient of $(X,\omega)$ by the hyperelliptic 
involution. Then the $\Sl(\R)$-orbit of $(\P^1, q)$ is dense in the 
stratum of meromorphic quadratic differentials having two simple zeores and six
simple poles. The stabilizer of the $\Sl(\R)$-action contains two
(non-commuting) pseudo-Anosov diffeomorphisms.
\end{Corollary}

One can reformulate the latter result in term of holomorphic quadratic
differentials on a genus $2$ surface in the following way. Let
$\QQQ_2$ be the moduli space of holomorphic quadratic differentials which
are {\it not} the global square of any Abelian differentials. Let $\pi
: Y \rightarrow \P^1$ be a double covering of $\P^1$ ramified
precisely over the six poles of $q$. Let $(Y,\tilde{q})$ be the
lift of $q$. Obviously  $(Y,\tilde{q}) \in \QQQ_2$. 

\begin{Corollary} 
\label{quadratic-diff2}
The $\Sl(\R)$-orbit of $(Y ,\tilde{q})$ is dense inside the whole
moduli space $\QQQ_2$. Moreover the stabilizer of the $\Sl(\R)$-action of
$(Y ,\tilde{q})$ contains two non-commuting pseudo-Anosov
diffeomorphisms.
\end{Corollary}

\begin{proof}
As we have seen, $(\P^1,q)$ belongs to the stratum $\QQQ(1, 1, -1^6)$ of
meromorphic quadratic differentials having two simple zeores and six 
simple poles. The stratum $\QQQ(1,1,-1^6)$ is isomorphic to the stratum 
$\QQQ(1,1,1,1)$, the principal stratum of holomorpic quadratic
differentials on genus $2$ surfaces (see Lanneau~\cite{La}). Moreover
the $\Sl(\R)$-action is equivariant with respect to this
isomorphism. Therefore Corollary~\ref{quadratic-diff1} shows that 
the closure of the Teichm\"uller disc $\Sl(\R)\cdot(Y ,\tilde{q})$ 
contains $\QQQ(1,1,1,1)$. This last stratum is dense inside $\QQQ_2$ 
which gives Corollary \ref{quadratic-diff2}.
\end{proof}

\begin{Acknowledgments}
We would like to thank the anonymous referee for a careful study 
and valuable remarks in order to clarify the paper. We are also grateful 
to Essen University where this project started.
This work was partially supported by the ANR Teichm\"uller ``projet blanc'' ANR-06-BLAN-0038.
\end{Acknowledgments}

\section{Background}\label{sec:background}

In order to establish notations and preparatory material, we
review basic notions concerning translation surfaces, affine 
diffeomorphisms, and moduli spaces. 

\subsection{Pseudo-Anosov diffeomorphisms and Veech groups}

A {\it translation surface} is a genus $g$ surface with a translation 
structure (i.e. an atlas such that all transition functions are translations). 
As usual, we
consider maximal atlases. These surfaces are precisely those given
by a Riemann surface $X$ and a holomorphic (non-zero) one-form
$\omega \in \Omega(X)$; see say~\cite{MT} for a general reference
on translation surfaces and holomorphic one-forms.

Let $(X,\omega)$ be a translation surface.  The stabilizer of $(X,\omega)$ under the 
$\SL(2,\R)$-action is called the {\it Veech group} of $(X,\omega)$ and is denoted 
by $\SL(X,\omega)$. A more intrinsic definition is the following.  An 
{\it affine diffeomorphism} is an orientation preserving homeomorphism of $X$ 
which is affine in the charts of $\omega$ and permutes the zeroes of $\omega$.  
The derivative (in the charts of $\omega$) of an affine diffeomorphism defines 
an element of $\SL(X,\omega)$. Conversely such element is the
derivative of an affine diffeomorphism.

A diffeomorphism $f$ is a {\em pseudo-Anosov diffeomorphism} if and only if the 
linear map $D f$ is
hyperbolic; that is $|\textrm{trace}(D f)| > 2$. In this case, $D
f$ has two real eigenvalues $\lambda^{-1} < 1 < \lambda$. The number $\lambda$ is 
 called the  {\it expansion factor} of the
pseudo-Anosov diffeomorphism~$f$.
  
\subsection{Connected components of the strata}
\label{strata}
The moduli space of Abelian differentials is stratified by the 
combinatorics of the zeroes. We denote by $\HHH(k_1, \dots, k_n)$ 
the stratum of Abelian differentials consisting of holomorphic one-forms with $n$ 
zeroes of multiplicities $(k_1,\dots,k_n)$. These strata are 
non-connected in general but each stratum has at most three connected components 
(see~\cite{KZ} for a complete classification). In particular 
the stratum with two zeroes of multiplicity $2$, $\HHH(2,2)$, has 
two connected components. The {\it hyperelliptic component}
$\HHH^{hyp}(2,2)$ contains precisely pairs $(X,\omega)$ of a
hyperelliptic 
surface $X$ and a one-form whose zeros are interchanged by the hyperelliptic
involution.
The other ({\it non-hyperelliptic}) component $\HHH^{odd}(2,2)$ 
is distinguished by an odd parity of the spin structure. 
There are two ways to compute the parity 
of the spin structure of a translation surface $X$. The first way is to use 
the Arf formula on a symplectic basis (see~\cite{KZ}). The second possibility 
applies if $X$ comes from a quadratic differential, i.e.\ if $X$ 
possesses an involution  such that the quotient produces a half-translation 
surface (see~\cite{La2}). \medskip

Let $\mathcal Q(1,1,-1^6)$ be the stratum of meromorphic 
quadratic differentials 
on the projective line with two simple zeroes and six simple poles. It is 
easy to see that this stratum is connected (see~\cite{KZ,La}). Taking 
the orientating double covering, one gets a local embedding
$$
\mathcal Q(1,1,-1^6) \rightarrow \mathcal H(2,2).
$$
We will denote by $\LLL$ the image in $\mathcal H(2,2)$ of the previous map. 
The construction of the Arnoux-Yoccoz surface $(X,\omega)$ is given below.
Here we record for later use:

\begin{Lemma}
The Arnoux-Yoccoz surface $(X,\omega)$ lies in $\LLL \subset \HHH^{odd}(2,2)$.
\end{Lemma}

\begin{proof}
Thanks to the decomposition of the Arnoux-Yoccoz surface $(X,\omega)$ 
into cylinders (see Section~\ref{firstdir}) it is easy to define 
an affine diffeomorphism of $(X,\omega)$ by a rotation of
$180$ degree around the center of $C'_1$ (see Figure~\ref{fig:adjusted}).
This diffeomorphism fixes $8$ points, but not the two zeros of
$\omega$. Hence $X$ is hyperelliptic and lies in $\HHH^{odd}(2,2)$.
\par

Let us recall the formula in Theorem~1.2. of~\cite{La2} (p.~516) in order 
to calculate the parity of the spin structure. If $X$ is a half translation surface 
belonging to the stratum $\QQQ(k_1,\dots,k_l)$ then the parity of the 
spin structure of $\hat X$, the orientating surface, is
$$
\left[\cfrac{|n_{+1}-n_{-1}|}{4}\right]\ \mod 2
$$
where $n_{\pm1}$ is the number of zeros of degrees $k_j = \pm 1 \mod  4$, and
where all the remaining zeros satisfy $k_r = 0 \mod  4$ 
(the square brackets denote the integer part). 
This shows that the parity of the spin structure of $(X_3,\omega_3)$ 
is then $\cfrac{1}{4}(6-2)=1 \mod 2$.
\end{proof}

Therefore the hyperelliptic locus $\LLL$ belongs to 
the odd part $\HHH^{odd}(2,2)$ of $\mathcal H(2,2)$. 
We recall that the complex dimension of $\mathcal Q(1,1,-1^6)$ is 
$6$ and the complex dimension of $\HHH(2,2)$ is $7$.

\begin{Lemma}
There exists a linear isomorphism between the stratum 
$\QQQ(1,1,1,1)$ and the stratum $\QQQ(1,1,-1^6)$.
\end{Lemma}

Here linear implies in particular that the $\Sl(\R)$-action commutes 
with this isomorphism.

\begin{proof}
We recall here the proof presented in~\cite{La}.
Let us consider a meromorphic quadratic differential 
$q$ on $\P^1$ having the singularity pattern
$(1,1,-1^6)$. Consider a ramified double covering $\pi$ over 
$\P^1$
having ramification points over the simple poles of $q$, and no
other ramification points. We obtain a genus $2$ 
hyperelliptic Riemann surface $X$ with a quadratic differential
$\tilde q = \pi^{*} q$ on it. It is easy to see that the induced quadratic differential 
has the singularity pattern $(1,1,1,1)$. Hence we get locally an
$\SL_2(\R)$-equivariant mapping
$$
\mathcal Q(1,1,-1^6) \to \mathcal Q(1,1,1,1).
$$
Since the dimensions coincide, i.e.\
$$\dim_{\C}\mathcal Q(1,1,-1^6)=2\cdot 0 + 8 - 2=6 =2\cdot 2 + 4 -2 
= \dim_{\C}\mathcal Q(1,1,1,1) $$
and since the geodesic flow acts ergodically on the strata
the image of the above map equals $\mathcal Q(1,1,1,1)$.
\end{proof}

\section{Completely periodic directions}

Let $(X,\omega)$ be a translation surface. A {\em cylinder} is a
topological cylinder embedded in $X$, isometric to a flat cylinder
$\R / w \Z \times (0,h)$. The boundary of a maximal cylinder is a 
union of a finite number of saddle connections.
 
A direction $\theta$ is {\em completely periodic} on $X$ if 
all the regular
geodesics in the direction $\theta$  are closed. This means that $X$
is the closure of a finite number of maximal cylinders in the
direction $\theta$. In a periodic direction, all the geodesics 
emanating from singularities are saddle connections.
 
Let $\theta$ be a completely periodic direction on a translation surface
$(X,\omega)$. A translation surface comes with a horizontal
and vertical direction, and we will henceforth assume that $\theta$
is different from the vertical direction. The saddle connections in the 
direction $\theta$ are labeled by $\gamma_1, \dots,\gamma_k$. The cylinders
are labeled by $\CCC_1, \dots, \CCC_p$ and $w_1,\dots,w_p$ will stand for 
the widths (or perimeters) of the cylinders. 

For each cylinder $\CCC_i$ one can encode the sequence of saddle
connections contained in the bottom of $\CCC_i$ and ordered in
the cyclic ordering of the boundary of $\CCC_i$ by a cyclic
permutation $\sigma^b_i$. We get an analogous definition if we replace
bottom by top. Therefore one gets two $p-$tuples of cyclic permutations 
$(\sigma^b_1 \dots \sigma^b_p)$ and $(\sigma^t_1 \dots \sigma^t_p)$. 
Note that these data 
define two permutations on $k$ elements $\pi_b = \sigma^b_1 \circ
\dots \circ \sigma^b_p$ and $\pi_t = \sigma^t_1 \circ \dots \circ
\sigma^t_p$. 

These data form the {\it combinatorics} $\comb$ of the 
direction $\theta$ on
the surface $X$. This notion is very close to the one of separatrix diagram
introduced by Kontsevich and Zorich (see~\cite{KZ}) 
but will be more convenient for our purposes.

To give a complete description of the surface $(X,\omega)$ in the
direction $\theta$ we also need continuous parameters: 
\begin{itemize}
\item the lengths of the saddle connections,
\item the {\it heights} of the cylinders with respect to the vertical direction,
\item the {\it twists} of the cylinders.
\end{itemize}
  
The only parameters which are non-trivial to define are the
twists. For that, one has first to fix a {\it marking} on the
combinatorics $\comb$, i.e.\ on each cycle of $\sigma^b_i$ and 
$\sigma^t_i$, we mark an arbitrary element, denoted by $m_b(i)$
(resp. $m_t(i)$). 
When representing the combinatorics by a table, we will underline
the marked elements.

\begin{figure}
\begin{center}
\psfrag{1}{$\scriptstyle 1$} \psfrag{2mt}{$\scriptstyle 2=m_t(1)$} 
\psfrag{2mb}{$\scriptstyle 2=m_b(2)$} \psfrag{3}{$\scriptstyle 3$} 
\psfrag{3m}{$\scriptstyle 3=m_b(1)$} \psfrag{4}{$\scriptstyle 4$} 
\psfrag{4m}{$\scriptstyle 4=m_t(2)$}\psfrag{5}{$\scriptstyle 5$} 

\psfrag{ci1}{$\scriptstyle \overrightarrow{w_1}$} 
\psfrag{ci2}{$\scriptstyle \overrightarrow{w_2}$} 
\psfrag{ccc1}{$\scriptstyle \CCC_1$} \psfrag{ccc2}{$\scriptstyle \CCC_2$}
\psfrag{w1}{$\scriptstyle \overrightarrow{h_1}$}
\psfrag{w2}{$\scriptstyle \overrightarrow{h_2}$}
\psfrag{v1}{$\scriptstyle \overrightarrow{v_1}$}
\psfrag{v2}{$\scriptstyle \overrightarrow{v_2}$}
\psfrag{t2}{$\scriptstyle \overrightarrow{t_2} = 
\overrightarrow{v_2}- \frac{(h_2)_y}{(h_1)_y}\overrightarrow{v_1}$}
\psfrag{r}{$\scriptstyle \frac{(h_2)_y}{(h_1)_y}\overrightarrow{v_1}$}
\psfrag{w1+v1}{$\scriptstyle \overrightarrow{h_1}+\overrightarrow{v_1}$}


\begin{minipage}[c]{.46\linewidth}
        \includegraphics[width=7cm]{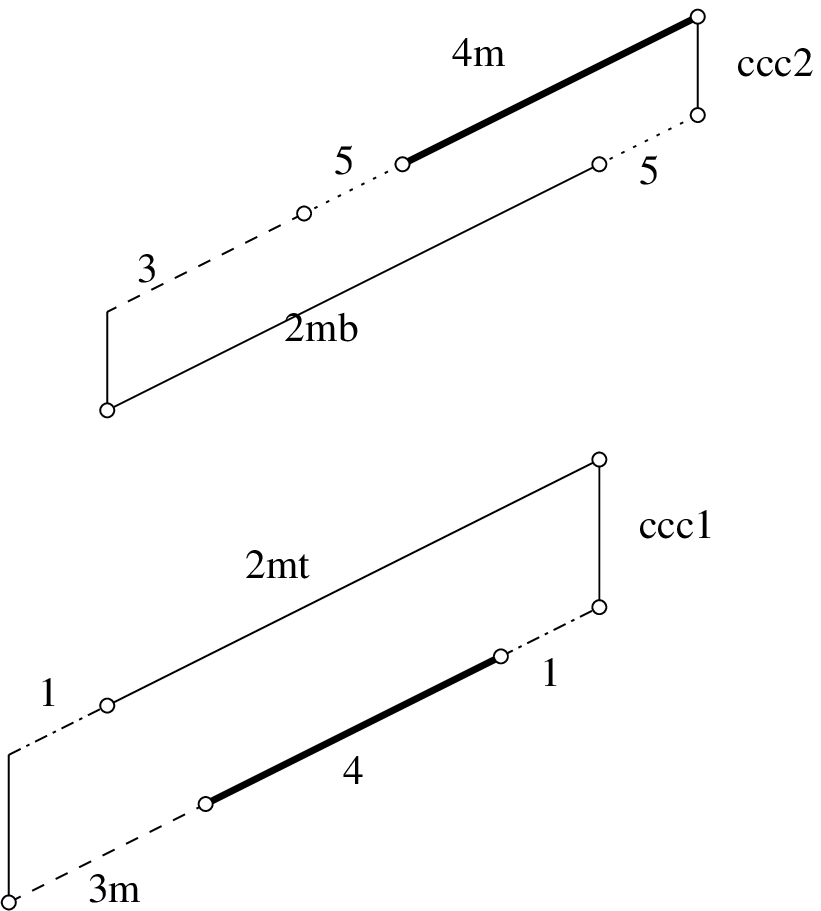}
   \end{minipage} \hfill
   \begin{minipage}[c]{.46\linewidth}
      \includegraphics[width=7cm]{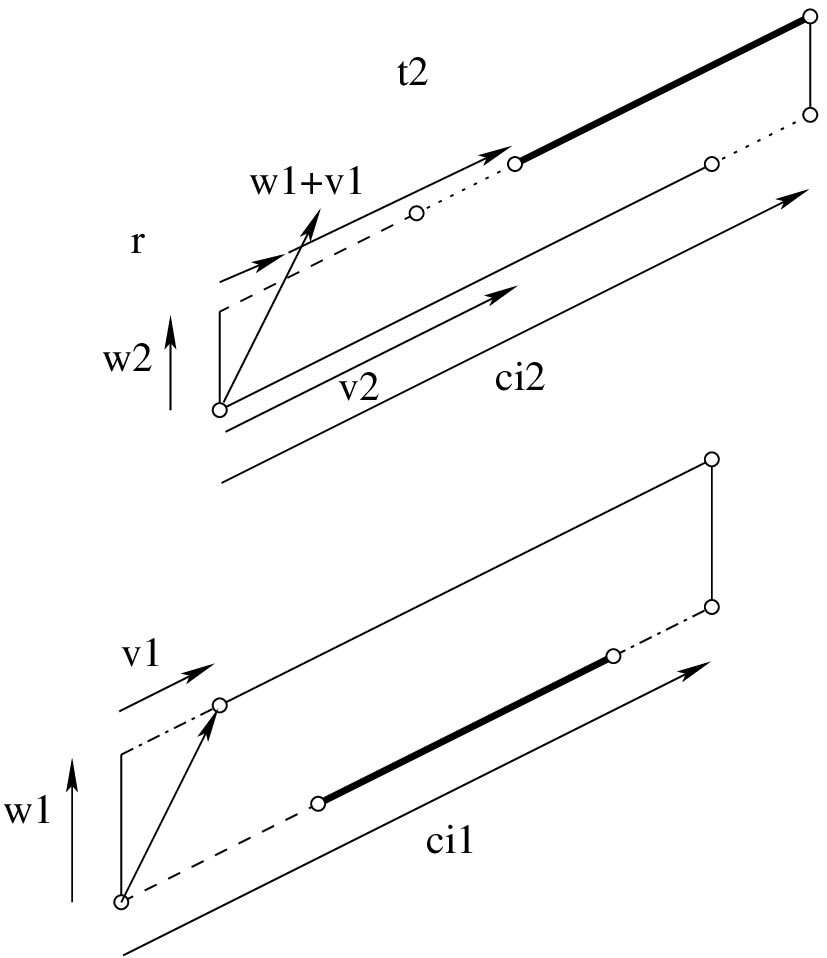}
   \end{minipage}

\end{center}
\caption{\label{fig:combinatorics}
Combinatorics and marking of a translation surface.
}
\end{figure}

On a translation surface with marked combinatorics in some
direction $\theta$ we will normalize the first twist to zero,
but define a twist vector in the cylinder $\CCC_1$. This 
is a saddle connection contained in $\CCC_1$ joining the origin of the saddle 
connection $\gamma_{m_b(1)}$ to the origin of the saddle connection
$\gamma_{m_t(1)}$. The corresponding vector may be decomposed into
its vertical component (denoted by $h_1$, where $(h_1)_y$ is the height of $\CCC_1$)
and its component in the direction
$\theta$, denoted by $v_1$ (see Figure~\ref{fig:combinatorics}).
The vector $v_1$ is well defined up to an additive constant $n w_1$ where 
$n\in \mathbb Z$. We normalize $v_1$ by requiring $|v_1| < |w_1|$ and
$v_1$ to be positive in the direction of $\theta$. 
\par
To have enough flexibility we will define twists for cylinders 
$\CCC_i$ ($i=2,\dots,p$) {\it with respect to the direction
$\theta^\perp = \theta^\perp(n_0)$}
given by $h_1 + v_1 +n_0 w_1$ for $n_0 \in \mathbb Z$. This is 
done in the following way. 
\par
Let $\CCC_i$ be a cylinder. Let $(h_i)_y$ be its height. The 
endpoint $P$ of the vertical vector 
$h_i = \left(\begin{smallmatrix}  0 \\ (h_i)_y \end{smallmatrix}\right)$ 
based at the origin of the saddle connection $\gamma_{m_b(i)}$ is located 
on the top of $\CCC_i$. Let $v_i$ be the vector joining $P$ 
to the origin of the saddle connection $\gamma_{m_t(i)}$ in the 
direction $\theta$. 
The vector $v_i$ is well-defined up an additive constant $n w_i$
where $n\in \mathbb Z$. The twist $t_i$ of $\CCC_i$ is defined to be the difference
$$
t_i = v_i - \frac{(h_i)_y}{(h_1)_y} (v_1 + n_0 w_1)
$$
(see Figure~\ref{fig:combinatorics}). The affine invariant will be the 
{\em normalized twist}, namely
$$
\frac{|t_i|}{|w_i|} \in [0,1[
$$
Therefore for each completely periodic direction $\theta$, each
marking $m$ on $\comb(\theta)$ and each $n_0\in \mathbb Z$, 
one gets the following quantities:
\begin{itemize}
\item $\overrightarrow{L}(\theta)\in \mathbb R^k$
\item $\overrightarrow{H}(\theta)\in \mathbb R^p$
\item $\comb(\theta)$
\item $\overrightarrow{T}(\theta,\marking,n_0)\in [0,1[^{p-1}$ (normalized twists)
\end{itemize}

Applying an appropriate linear transformation $M$ in $\textrm{GL}_2(\R)$ to the surface $(X, \omega)$ we can make all cylinders, in the directions $\theta$, horizontal. We send the direction $\theta^\perp$ to the vertical direction. Furthermore we normalize the cylinder $\CCC_1$ of direction $\theta$ to a unit square. This is the {\it normalizing  matrix} associated to the pair $(\theta, \theta^\perp)$.

Note that the conditions above define completely the linear transformation $M$. Let $h$ be the affine Dehn twist along the cylinder $\CCC_1$. If we change the direction $\theta^\perp$ by $h^{n_0}$ we then change the matrix by $M h^{-n_0}.$

 \begin{Remark}
\label{rk:combinatorics}
In Figure~\ref{fig:combinatorics} one has, with previous notations,
$\sigma^t_1 = (1\ 2), \sigma^t_2 = (3\ 5\ 4)$ and $\sigma^b_1=(3\ 4\ 1),\
\sigma^b_2=(2\ 5)$. Or equivalently: 
$\comb =
\left( \left(\begin{smallmatrix} 1&2 \\
  3&4&1\end{smallmatrix}\right),\ 
\left(\begin{smallmatrix} 3&5&4 \\ 2&5 \end{smallmatrix}\right) \right)$
  and a marking is 
$(\comb,\marking) =
\left( \left(\begin{smallmatrix} 1&\underline{2} \\
  \underline{3}&4&1\end{smallmatrix}\right),\ 
\left(\begin{smallmatrix} 3&5&\underline{4} \\
  \underline{2}&5\end{smallmatrix}\right) \right)$.
\end{Remark}

For a vector $u\in \mathbb R^n$ and a permutation $\pi$ on $n$ symbols, 
we will use the following convention 
$\pi \left( (u_1,\dots,u_n) \right) = (u_{\pi(1)},\dots,u_{\pi(n)})$.

We then have the obvious 

\begin{Lemma}
\label{lm:invariant:1}
Let $f\ :\ X\longrightarrow X$ be an affine diffeomorphism. Let
$\theta \in \mathbb S^1$ be a completely periodic direction and
$\theta':=f(\theta)$. Then $f$ induces a bijection between the
collections $\{\gamma_1,\dots,\gamma_k\}$ and 
$\{\gamma'_1,\dots,\gamma'_k\}$ of saddle connections 
(and thus a permutation $\pi^{sc}_f$ on $k$ elements) and a bijection
between the collections $\{\CCC_1,\dots,\CCC_p\}$ and $\{\CCC'_1,\dots,\CCC'_p\}$
of cylinders (and thus a permutation $\pi^{cy}_f$ on
$p$ elements). Moreover one has:
$$
\overrightarrow{L}(\theta) = \pi^{sc}_f \left( \overrightarrow{L}(\theta')
\right) \ \in \mathbb{RP}(k)
$$
and
$$
\overrightarrow{H}(\theta) = \pi^{cy}_f \left( \overrightarrow{H}(\theta')
\right) \ \in \mathbb{RP}(p).
$$
\end{Lemma}

\begin{Lemma}
With the same assumptions as in Lemma~\ref{lm:invariant:1}, let us 
choose a marking $\marking$ on $\comb = \comb(\theta)$. Then 
$\pi^{sc}_f$ induces a marking $\marking'$ on $\comb' = \comb(\theta')$.
Moreover there exists $n'_0\in \mathbb Z$ such that the
normalized twists $\overrightarrow{T}(\theta,\marking,0)$ and 
$\pi^{sc}_f \left( \overrightarrow{T}(\theta',\marking',n'_0)
\right)$ are the same.
\end{Lemma}

In fact, given a combinatorics in the direction $\theta$, the lengths of the
saddle connections, the width of the cylinders (with respect to the
vertical direction) and the twists characterize the surface
$(X,\omega)$ in the moduli space. Namely one has:

\begin{Theorem}
\label{theo:invariants}
Let $X$ be a flat surface. Let $\theta,\theta' \in \mathbb S^1$ be two
completely periodic directions. Let us choose a marking $\marking$ on
$\comb$. Let us also assume there exist two permutations 
$\pi_1$ on $k$ elements and $\pi_2$ on $p$ elements such that
\begin{enumerate}

\item $\overrightarrow{L} = \pi_1 \left( \overrightarrow{L'} \right) \
\in \mathbb{RP}(k)$ and $\overrightarrow{H} = \pi_2 \left(
\overrightarrow{H'} \right) \ \in \mathbb{RP}(p)$,

\item $\comb$ and $\comb'$ are isomorphic via $(\pi_1,\pi_2)$ and 
there exists $n'_0\in\mathbb Z$ such that the normalized twists
$\overrightarrow{T}(\theta,\marking,0)$ and $\pi_1 \left(
\overrightarrow{T'}(\theta',\marking',n_0') \right)$ are the same.

\end{enumerate}

Then there exists an affine diffeomorphism $f \in \textrm{Aff}(X,\omega)$. 
Moreover $f(\theta)=\theta'$, $f(\theta'^\perp)=\theta'^\perp$ 
and $\pi^{sc}_f = \pi_1,\ \pi^{cy}_f = \pi_2$.

Moreover if $M$ and $M'$ are the normalizing matrices associated with the 
pairs $(\theta, \theta^\perp)$ and $(\theta', \theta'^\perp)$, then $D f = M'^{-1}M.$
\end{Theorem}


\section{Application to the Arnoux-Yoccoz flat surface}
\label{sec:calculus} 
\subsection{Construction}(after P.~Arnoux)

Here we present the idea of the construction of the Arnoux-Yoccoz 
flat surface. What appears to be ad hoc here was constructed
originally by a self-similar interval exchange transformation
whose canonical suspension is a quadratic differential on $\P^1$. 
One can find all details in the original 
paper~\cite{Arnoux:Yoccoz}. See also~\cite{Ar2}, pp.~ 496--498. 
\par
\begin{figure}
\begin{center}
\psfrag{A}{$\scriptstyle A$}
\psfrag{B}{$\scriptstyle B$} \psfrag{Bp}{$\scriptstyle B'$}
\psfrag{C}{$\scriptstyle C$} \psfrag{Cp}{$\scriptstyle C'$}
\psfrag{D}{$\scriptstyle D$} \psfrag{Dp}{$\scriptstyle D'$}

\psfrag{1}{$\scriptstyle 1$}  \psfrag{2}{$\scriptstyle 2$} 
\psfrag{3}{$\scriptstyle 3$}  \psfrag{4}{$\scriptstyle 4$} 
\psfrag{5}{$\scriptstyle 5$}  \psfrag{6}{$\scriptstyle 6$} 
\psfrag{7}{$\scriptstyle 7$}  

\includegraphics[width=8.5cm]{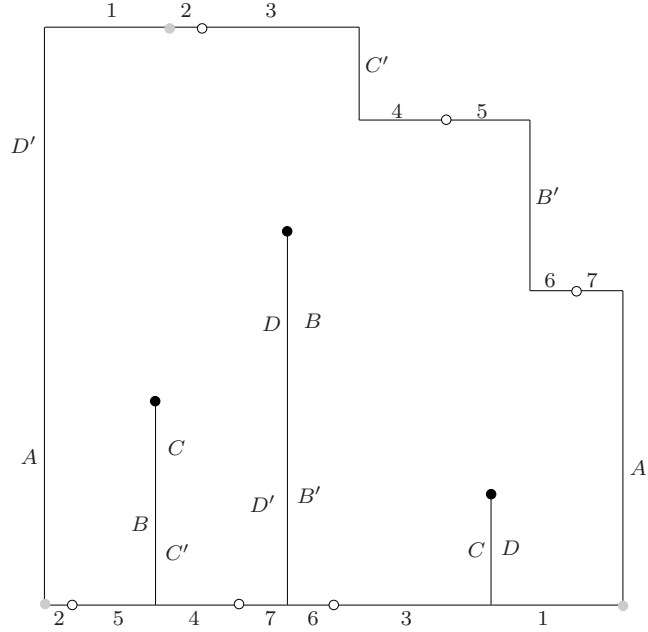}
\end{center}
\caption{\label{fig:gluings}
The Arnoux-Yoccoz flat surface. Identifications of the boundaries 
are made with respect to the labels. The black and white bullets correspond 
to the two singularities. The grey bullets correspond to a marked point.
}
\end{figure}
\par
Let us consider three rectangles in $\mathbb R^2$ glued along the
horizontal segment $[0,2]$ with the following parameters.
The bases have length $2\alpha,\ 2\alpha^2$ and $2\alpha^3$ and 
the heights are respectively $2,\ 2(\alpha+\alpha^2)$ and $2\alpha$ 
(here $\alpha$ is the positive real root of polynomial $X^3+X^2+X-1$). 
On the union of these three rectangles, we make  small vertical cuts starting 
on the base with  specified lengths. The cuts start at the points
$\alpha - \alpha^3,\ 
\alpha + \alpha^2,\ 1+ \alpha$ respectively on the base and
they have heights $\alpha +
\alpha^3,\ 1 + \alpha^2,\ \alpha^2 + \alpha^4$ respectively 
(see Figure~\ref{fig:gluings}). 
\par
We then identify horizontal boundaries with the following rules: we 
identify the point $(x,y)$ on the top boundary with $(f(x),0)$ for any 
$x\in [0,2]$, where $f: [0,2[ \to [0,2[$ is an interval exchange 
transformation. Following \cite{Ar2}~page 489, this interval exchange 
transformation has permutation $\pi=\left(
\begin{smallmatrix} 1&2&3&4&5&6&7\\ 2&5&4&7&6&3&1\end{smallmatrix}\right)$
and the interval labels are indicated in Figure~\ref{fig:gluings}.
\par
In this manner we obtain a flat surface with vertical boundary 
components. To get a closed flat surface we have to identify the 
vertical boundaries. This is done according to the rules presented in 
Figure~\ref{fig:gluings}. The parameters are chosen in such a way that the 
gluing is by (Euclidean) isometry.
\par
The flat surface $(X,\omega)$ obtained in this way has
two conical singularities of
total angle $6\pi$, so $(X,\omega) \in \mathcal H(2,2)$.
\par 
We review the construction of the pseudo-Anosov diffeomorphism
due to Arnoux-Yoccoz. Let us
define a new parametrization of $X$. We will do that in two
steps. First in the chart of $X$, let us 
take the two small rectangles and glue them above the big rectangle by
an isometry (in this operation, we have to cut the medium rectangle
into two parts). Then we cut the figure along the vertical line
$x=\alpha+\alpha^4$ and permute the two obtained pieces. In
this way we get a new parametrization of our surface which is 
exactly the same except that the horizontal coordinates are
multiplied by $\alpha$ and the vertical coordinates are multiplied by
$\alpha^{-1}$. Therefore we can define an affine diffeomorphism~$\Phi$ on 
$(X,\omega)$:
$$
\Phi(x,y) \textrm{ in the first parametrization} = (\alpha
x,\alpha^{-1}y) \textrm{ in the second parametrization.}
$$
The derivative of $\Phi$ is
$$
\left(\begin{array}{cc}
\alpha & 0 \\
0 & \alpha^{-1}
\end{array}\right) = 
\left(\begin{array}{cc}
\lambda^{-1} & 0 \\
0 & \lambda
\end{array}\right),\ \frac1{\alpha}=\lambda=\lambda(\Phi) > 1
$$

\subsection{Proof of Theorem \ref{pseudo-Anosov}}

There are two completely directions on the Arnoux-Yoccoz surface of slopes $\theta = 1- \alpha^2$ and $\theta' = 3 + \alpha^2$ respectively. Let $\theta^\perp = 1+\alpha^2$ and $\theta'^\perp = \frac1{169}(367+252\alpha+175\alpha^2)$. We check that the 
hypothesis of Theorem \ref{theo:invariants} for the pairs $(\theta, \theta^\perp)$, $(\theta', \theta'^\perp)$ are fulfilled. The normalizing matrices associated to the pairs $(\theta, \theta^\perp)$, $(\theta', \theta'^\perp)$ are:

$$M = 
\frac1{4}\left(\begin{smallmatrix} 
1+\alpha^2  &   & -4 -4\alpha-2\alpha^2\\ 
-1 &  & 6+5\alpha +3\alpha^2
\end{smallmatrix}\right) \qquad
M' = 
\frac1{4}\left(\begin{smallmatrix} 
-49-42\alpha-27\alpha^2  &   & -66 -56\alpha-36\alpha^2\\ 
14+14\alpha+9\alpha^2 &  & 20+17\alpha +11\alpha^2
\end{smallmatrix}\right).$$
Therefore there exists a pseudo-Anosov diffeomorphism $\tilde{\Phi}$ on the Arnoux-Yoccoz flat surface whose differential is
$$
D \tilde \Phi = M'^{-1}M = \left(
\begin{smallmatrix} 23+18\alpha+12\alpha^2 & &-29-24\alpha-16\alpha^2 \\ \\
74+62\alpha+40\alpha^2 & & -95-80\alpha-52\alpha^2 
\end{smallmatrix}
\right).
$$
\medskip
Since the trace of this matrix is greater than $2$, the diffeomorphism $\tilde
\Phi$ is a pseudo-Anosov diffeomorphism. The minimal polynomial over $\Q$ of the expansion factor of $\tilde \Phi$ is $$X^6+114X^5-409X^4+604X^3-409X^2+114X+1.$$
Moreover neither the horizontal direction
nor the vertical direction is an eigenvector of $D \tilde \Phi$.
The eigenvectors of the derivative correspond to the stable
and unstable foliation. Since the stable and unstable foliation
of $\Phi$ are by construction the horizontal and vertical one
we conclude that both the stable and unstable foliations of $\Phi$ 
and $\tilde \Phi$ are different. Since two pseudo-Anosov
diffeomorphisms have a power in common if and only if their 
stable and unstable foliations coincide, this conclude the
proof of Theorem \ref{pseudo-Anosov}. The details of the calculations will be given in the Appendix.

\section{Closure of the disc}

The set of translation surfaces $(X,\omega) \in \LLL$ with fixed
area $\alpha$ (with respect to $\omega$) form a real hypersurface
$\LLL_\alpha$ in $\LLL$.
\par
This section is devoted to prove

\begin{Theorem} \label{closure}
The Teichm\"uller disc stabilized by the Arnoux-Yoccoz
pseudo-Anosov diffeomorphism is dense inside $\LLL_\alpha$.
\end{Theorem}

Equivalently, we claim that the $\GL^+_2(\R)$-orbit of the Arnoux-Yoccoz
surface is dense inside $\LLL$. We will switch between
both statements in the proof: Using $\GL^+_2(\R)$ is more
convenient for normalization while we need an ergodicity
argument that holds for the $\SL_2(\R)$-action on $\LLL_\alpha$.
\par
We sketch the strategy of proof. Our proof follows the ideas introduced by McMullen in \cite{Mc3}. First, by topological
considerations we find a direction on the
Arnoux-Yoccoz surface with a set of
homologous saddle connections that allows one to decompose the 
surface into tori and cylinders. The parabolic subgroup of 
$\SL_2(\R)$ that  fixes the direction of these saddle connections acts on the 
set of pairs of lattices that define the tori. 
Second we apply Ratner's theorem and the explicit geometry
of the surface to see that the orbit closure in the homogeneous
space of pairs of lattices is as big as one could hope for.
Hence the closure of the disc contains at least
the submanifold of $\LLL_\alpha$ defined by fixing the areas of
the decomposition pieces. We will give details of this
using period coordinates on $\LLL$. The third step consists
of splitting the Arnoux-Yoccoz surface in another direction
to remove this area constraint.
\par
Let $(X,\omega)$ be the Arnoux-Yoccoz surface. Besides the completely periodic 
directions (for example the one 
with slope $\theta=1-\alpha^2$
used in the previous section) there are directions $\Xi$ on $(X,\omega)$ 
with the following property. There are $4$ homologous saddle connections 
$\beta_i$ such that $\beta_2$ and $\beta_3$ (resp.\ $\beta_4$ and $\beta_1$)
bound a cylinder $C_1$ (resp.\ $C_2)$. The hyperelliptic involution
interchanges $C_1$ and $C_2$. The complement of $X \setminus \{C_1 \cup C_2\}$
consists of two components $T_1$ and $T_2$. If we identify the cuts
lines $\beta_1$ and $\beta_2$ on $T_1$ and $\beta_3$ 
and $\beta_4$ on  $T_2$ we obtain two tori, which  also denoted by the 
same letter. In fact $T_1$ and $T_2$ are irrationally foliated as
we will see below. 
\par
\begin{figure}
\begin{center}

\includegraphics[width=8.5cm]{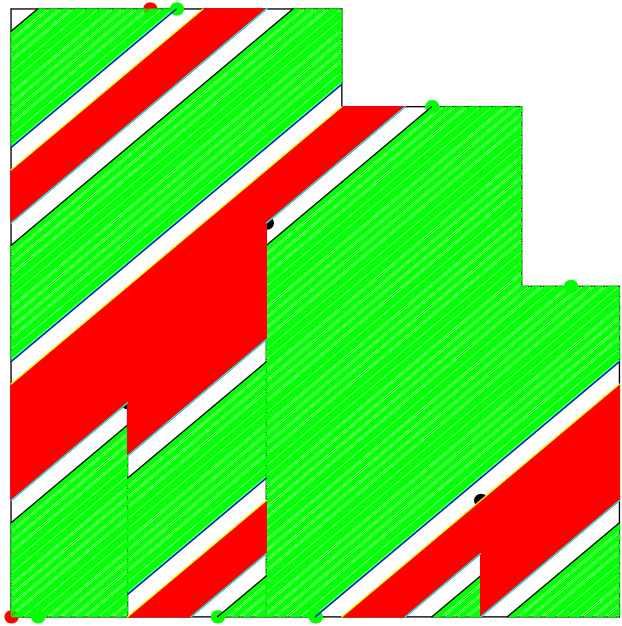}
\end{center}
\caption{\label{fig:direction:2T2C}
Decomposition of the Arnoux-Yoccoz surface: shaded regions
represent a 2T2C-direction, black lines the 3C-direction
studied in Section~\ref{firstdir}.
}
\end{figure}

Such a direction is given by the slope $\Xi = \alpha+\alpha^2$ (see
Figure~\ref{fig:direction:2T2C}). 

\begin{figure}
\begin{center}

\includegraphics[width=11cm]{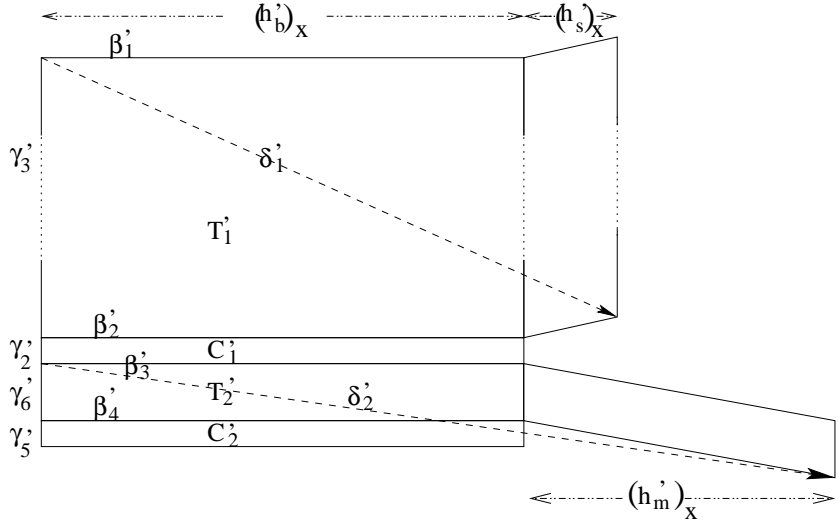}
\end{center}
\caption{\label{fig:adjusted}
The adjusted Arnoux-Yoccoz flat surface.
}
\end{figure}

\begin{Convention}Abusing notations, we will use the same letter for a saddle connection and its affine holonomy vector. 
\end{Convention}
We apply 
a matrix in $\GL^+_2(\R)$ such that $\Xi$ becomes horizontal
and $\theta$ vertical and moreover such that $\gamma_2$ becomes
$\gamma'_2 = \left(\begin{smallmatrix} 0 \\ 1 
\end{smallmatrix}\right)$ and $\beta_i$ becomes 
$\beta'_i = \left(\begin{smallmatrix} 1 \\ 0 
\end{smallmatrix}\right)$. We can draw this surface in the 
$\GL^+_2(\R)$-orbit of $(X,\omega)$ as indicated in Figure \ref{fig:adjusted}.
Note that  the primes in this section do not have the same meaning 
as the primes in Section \ref{seconddir}.  Note moreover that the figure
does not display scales in the vertical and horizontal direction
correctly.
\par
\par
We  refer to this surface as the adjusted Arnoux-Yoccoz surface
$(X',\omega')$.
Moreover we call a splitting of a translation surface in $\LLL$
with the same topology and dynamics as the horizontal one a 
{\em $2T2C$-splitting}. 
\par

\subsection{Applying Ratner's theorem}

Let $G$ be a Lie group, $\Gamma$ a lattice in $G$ and $U \subset G$  
a $1$-parameter subgroup generated by unipotent elements. The group $U$ acts 
on the left on $\XXX = G/\Gamma$.  Ratner's theorem  (\cite{Ra}) states 
that for any $x \in \XXX$ the closure $\ol{U \cdot x}$
is an orbit $H \cdot x \in \XXX$, where $H$ is a unimodular subgroup
depending on $x$ with the property that $x \Gamma x^{-1} \cap H$ 
is a lattice in $H$. \medskip

Let $\Lambda_i$, $i=1,2,3$ be lattices in $\R^2$ normalized
such that $${\rm area}(\R^2/\Lambda_i)=1.$$
Triples of normalized lattices are parameterized by the homogeneous space
$\XXX = G^3/ \Gamma^3$ where $G= \SL_2(\R)$ and $\Gamma = \SL_2(\Z)$.
We denote the projections onto the factors resp.\ pairs of
factors by ${\rm pr}_i$ resp.\ ${\rm pr}_{ij}$. 
Some more notation: Let $\fracg$ be the
Lie algebra of $G$ and $\fracn$ (resp.\ $\fracu$, resp.\ $\fraca$) be 
the Lie algebra of the unipotent upper triangular matrices $N$ (resp.\
unipotent lower triangular matrices $U$, resp.\ diagonal matrices $A$).
Choose standard generators 
$$n =  \left(\begin{array}{cc} 0 & 1 \\ 0 & 0 \\ \end{array}
\right) \in \fracn, 
\quad a = \left(\begin{array}{cc} 1 & 0 \\ 0 & -1 \\ \end{array}
\right)  \in \fraca, 
\quad u = \left(\begin{array}{cc} 0 & 0 \\ 1 & 0 \\ \end{array}
\right) \in \fracu $$
such that $[n,a]=-2n$, $[n,u]=a$ and $[a,u]=-2u$. 

\par
We will apply this to the case where $U$ is the group diagonal
embedding $N_\Delta$ of the unipotent upper triangular matrices $N$.
We remark that, in genus 2, McMullen  studies
the case of one lattice and the action of $N  \cap \SL_2(\Z)$
or the case of two lattices and the action of $N_\Delta$
(see \cite{Mc3} Theorems~2.3 and 2.6 of loc.\ cit.). We remark that we 
could also proceed with Shah's version of Ratner's theorem for
cyclic groups (\cite{Sh}) and the action of $\N_\Delta \cap \SL_2(\Z)$
on pairs of lattices, fixing $C_i''$.
We will not list all the possible closures of $N_\Delta$-actions but 
restrict to what we actually need.
\par
\begin{Lemma} \label{3summands}
Suppose that $\Lambda_3$ is the standard lattice and 
$C :=  \ol{N_\Delta \cdot (\Lambda_1,\Lambda_2,
\Lambda_3)}$ projects
to the whole space $G^2/\Gamma^2$ via ${\rm pr}_{12}$.
Then $C = (G\times G \times N) \cdot (\Lambda_1,\Lambda_2, \Lambda_3)$. 
\end{Lemma}
\par
\begin{proof} By Ratner's theorem $C=H \cdot (\Lambda_1,\Lambda_2, 
\Lambda_3)$ for some $H$.
By the hypothesis on $\Lambda_3$ the Lie
algebra $\frach$ of $H$ is contained in 
$$(\fracn_1 \oplus \fraca_1 \oplus
\fracu_1) \oplus  (\fracn_2 \oplus \fraca_2 \oplus
\fracu_2) \oplus \fracn_3.$$ 
Either this is an equality and we are done
or $\frach$ is given by one equation
$\sum_{i=1}^7 \alpha_i x_i=0$, in the basis $\{e_1, \dots, e_7\}$ where $e_i$ are the standard generators of 
the summands defined above. For $i \neq 7$ let $b_i=e_i$ if $\alpha_i=0$
and $b_i = \alpha_7/\alpha_i e_i -  e_7$ otherwise. We have
$b_i \in \frach$ in both cases. One checks that
$[b_1,b_3]$ is a non-zero multiple of $e_2$, that
$[b_1,b_2]$ is a non-zero multiple of $e_1$ and that
$[b_3,b_2]$ is a non-zero multiple of $e_3$. Continuing
like this $e_i \in \frach$ for $i \neq 7$. Since $\frach$ contains
the diagonal we also have $e_7 \in \frach$ and we are done.
\end{proof}
\par
\begin{Corollary} \label{incomm2}
Let $\Lambda_3$ be the standard lattice.
Suppose that neither of the lattices $\Lambda_1$ and $\Lambda_2$ 
contains a horizontal vector and suppose there does not exist 
an element $M_t \in N$ such that $\Lambda_1$ and $M_t \cdot
\Lambda_2$ are commensurable. Then
$$ C:= \ol{N_\Delta \cdot (\Lambda_1,\Lambda_2, \Lambda_3)} = (G \times G
\times N)
\cdot (\Lambda_1,\Lambda_2, \Lambda_3) $$
\end{Corollary}
\par
\begin{proof} The hypothesis on $\Lambda_1$ and $\Lambda_2$ are
just the one needed to apply Theorem~2.6 in \cite{Mc3}. Its conclusion
is the condition on ${\rm pr}_{12}$ needed to apply Lemma~\ref{3summands}.
\end{proof}
\par
\smallskip
Let $C_1'$ and $T_i'$, $i=1,2$, be the components obtained by
splitting the adjusted Arnoux-Yoccoz surface $(X',\omega')$.
We denote by the same symbols the tori obtained by gluing the
slits. Let $\Lambda'_i$, $i=1,2$ and $\Lambda'_C$ be 
defined by $\R^2/\Lambda'_i \cong T'_i$ and 
$\Lambda'_3 := \Lambda'_C$ be defined by $\R^2/\Lambda'_C
\cong C'_1$. Finally we apply homotheties to $\Lambda'_i$ 
in order to obtain lattices $\widetilde{\Lambda_i}$ with area one.
We will generally denote with a tilde lattices that 
are {\em area-normalized}.
\par
\begin{Lemma} \label{Ratnerapplies2}
The area-normalized lattices $\widetilde{\Lambda_i}$ obtained by splitting 
the adjusted Arnoux-Yoccoz 
surface $(X',\omega')$ in the horizontal direction, i.e.\ along 
the saddle connections $\beta'_i$, 
satisfy the conditions of Corollary \ref{incomm2}. \\
More generally, let $\Lambda_1, \Lambda_2$ be the lattices 
$$
\Lambda_i = \langle \left(\begin{smallmatrix} 0 \\ a_i 
\end{smallmatrix}\right), \left(\begin{smallmatrix} b_i \\ c_i 
\end{smallmatrix}\right) \rangle .
$$
Let $A_i=\sqrt{{\rm area}(\R^2/\Lambda_i)}$ with $i=1,2$. Assume that $a_i, b_i, c_i$ belong to $\Q(\alpha)$ for $i=1,2$. If
$\cfrac{c_i}{a_i}\not \in \Q$ and $\Q(A_1)\neq \Q(A_2)$
then the {\it area-normalized}
lattices $\widetilde{\Lambda_1}$ and 
$\widetilde{\Lambda_2}$ 
satisfy the conditions of Corollary~\ref{incomm2}.
\end{Lemma}
\par
\begin{proof} We prove the general statement first. Nonexistence
of horizontal vectors in $\widetilde{\Lambda_i}$ (and thus in $\Lambda_i$) 
is equivalent to the first condition on
the lattices $\Lambda_i$. Now let us prove 
that the area-normalized lattices $\widetilde{\Lambda_1}$ and
$M_t\cdot\widetilde{\Lambda_2}$ are not commensurable for any $t\in \R$.
By contradiction let us assume there exist integers $n,m$ and $p$ such that
$$
n \left(\begin{matrix} 0 \\ \frac{a_1}{A_1} 
\end{matrix}\right) = m M_t \left(\begin{matrix} 0 \\ \frac{a_2}{A_2} 
\end{matrix}\right)
+p M_t \left(\begin{matrix} \frac{b_2}{A_2} \\ \frac{c_2}{A_2}
\end{matrix}\right).
$$
This yields
$$
n \left(\begin{matrix} 0 \\ \frac{a_1}{A_1} 
\end{matrix}\right) = m \left(\begin{matrix} \frac{t a_2}{A_2} \\ \frac{a_2}{A_2} \end{matrix}\right)
+p  \left(\begin{matrix} \frac{b_2+ t c_2}{A_2} \\ \frac{c_2}{A_2}
\end{matrix}\right).
$$
Taking the second coordinate one gets:
$$
n \frac{a_1}{A_1} = m  \frac{ a_2}{A_2} + p \frac{c_2}{A_2}
$$
implying $A_1 \in \Q(A_2)$ which is a contradiction. \bigskip

Now let us prove that the area-normalized lattices $\widetilde{\Lambda_i}$ 
obtained by splitting the adjusted Arnoux-Yoccoz 
surface satisfy the conditions of Corollary \ref{incomm2}. For that 
we will use the general statement we have just proved. \medskip

\begin{Claim}
The lattices $\Lambda'_1$ and $\Lambda'_2$ (and thus 
$\widetilde{\Lambda_1}$ and $\widetilde{\Lambda_2}$) have no 
horizontal vectors. 
\end{Claim}

\begin{proof}[Proof of Claim]
This exactly means that the horizontal flows on
$T'_1$ and $T'_2$ are irrational.  
We explain the method for $T'_1$. Flowing in the horizontal direction, from $P$ the intersection of
$\gamma'_3$ and $\beta'_1$ (see Figure~\ref{fig:adjusted}),
we cross again $\gamma'_3$ at a point $Q$. The flow in the horizontal
direction is irrational if and only if the first return map on the
circle $\gamma'_3$ is an irrational rotation. We prove that the
normalized twist $\frac{PQ}{\gamma'_3} \in \Q(\alpha)$ is irrational.  
A direct computation shows that the value of the normalized twist for
$T'_1$ is $\frac{1}{2} - \frac{\alpha^2}{2}$ and for $T'_2$ is  
$\frac{10}{11} - \frac{5}{11 \alpha} + \frac{3}{11 \alpha^2}$. Both
numbers are irrational which proves the claim. \medskip
\end{proof}

Now let us verify the second hypothesis. We have to prove that the numbers 
$\sqrt{{\rm area}(T'_1)}$ and $\sqrt{{\rm area}(T'_2)}$
generate two distinct quadratic extensions of $\Q(\alpha)$.
According to Figure~\ref{fig:adjusted} one has 
\begin{equation} \label{Tareas}
{\rm area}(T'_1) = \gamma'_3 \wedge \beta'_1 + \gamma'_3 \wedge h'_s
\quad \text{and} \quad
{\rm area}(T'_2) = \gamma'_6 \wedge \beta'_4 + \gamma'_6 \wedge
h'_m,
\end{equation}
where $\wedge$ is the cross product.
The adjusted Arnoux-Yoccoz surface is obtained from the original
one by the linear map which sends the vectors $\gamma_2$ to
$\left(\begin{smallmatrix} 0 \\ 1  
\end{smallmatrix}\right)$  and $\beta_2$ to
$\left(\begin{smallmatrix} 1 \\ 0  
\end{smallmatrix}\right)$. This is a calculation  in the field
$\Q(\alpha)$. Since 
$$
\frac{{\rm area}(T'_1)}{{\rm area}(T_1)} =
\frac{{\rm area}(T'_2)}{{\rm area}(T_2)} \in \Q(\alpha),
$$
we may as well prove the same assertion for $\sqrt{{\rm area}(T_1)}$ 
and $\sqrt{{\rm area}(T_2)}$ using Formula~\ref{Tareas} without the
primes and the values from table in Section~\ref{firstdir}. This leads to:
$$
{\rm area}(T_1)= 4 \alpha \qquad \textrm{and} \qquad 
{\rm area}(T_2) =  -8 + 12 \alpha + 8 \alpha^2.
$$
and the proof of Lemma \ref{Ratnerapplies2} is completed
by the following claim.
\end{proof}

\begin{Claim} \label{extension}
The numbers $\sqrt{{\rm area}(T_1)}$ and $\sqrt{{\rm area}(T_2)}$
generate two different quadratic extensions of $\Q(\alpha)$. 
\end{Claim}

\begin{proof}[Proof of the Claim]
We first prove that $x:=\sqrt{{\rm area}(T_1)}$ does not 
belong to $\Q(\alpha)$. The polynomial $Q(X) = X^6 + X^4 + X^2 - 1$
annihilates $x$ by definition of $\alpha$ and it is irreducible 
over $\Q$. Indeed 
$Q(X+1) = 2 + 12X + 22X^2 + 24 X^3 + 16 X^4 + 6 X^5 + X^6$ satisfies
Eisenstein's criterion with respect to the prime number $2$. 
\par
Second, assume by contradiction, that $y := \sqrt{-2 +3\alpha+2\alpha^2}$ 
belongs to $\Q(x)$. Then, there
exist $a,b \in \Q(\alpha)$ such that $y = a + bx$. Taking the square of
this equation we get $a^2 + b^2 \alpha + 2ab x = -2
+3\alpha+2\alpha^2$. In the basis $\{1,x\}$ this leads to the two following 
equations.
$$
\left\{\begin{array}{lcl}
ab & = & 0 \\
a^2 + b^2 \alpha & = & -2 +3\alpha+2\alpha^2
\end{array}\right.
$$
If $b = 0$, we get $a^2  =  -2 +3\alpha+2\alpha^2$. We have 
already proved that $-2 +3\alpha+2\alpha^2$ is not a square in 
$\Q(\alpha)$ which is a contradiction with $a \in \Q(\alpha)$. \\
If $a = 0$, we get $b^2 = 1 - 2 \alpha^2$. A straightforward 
computation shows that $b$ is a root of the polynomial 
$S(X)  = X^6 - 5X^4 + 19X^2 - 7$. This polynomial is irreducible over 
$\Q$. Consequently $\Q(b)$ has 
degree $6$ over $\Q$ which is a contradiction to $b \in \Q(\alpha)$.
The claim is proven. 
\end{proof}

\subsection{Splitting in different directions}

We want to conclude the proof of Theorem \ref{closure}.
By the following observation it suffices to show that
the closure of the Arnoux-Yoccoz disc contains a set of positive measure of $\LLL_\alpha$
or equivalently that the $\GL^+_2(\R)$-orbit of $(X',\omega')$
contains a set of positive measure of $\LLL$.
\par
\begin{Lemma} \label{ergodo}
Let $\LLL_\alpha$ be the hypersurface 
of fixed area $\alpha$ in the hyperelliptic locus $\LLL \subset \HHH(2,2)$.
The action of $\SL_2(\R)$ on  $\LLL_\alpha$
is ergodic.
\end{Lemma}
\par
\begin{proof} Let $\mathcal Q(1,1,-1^6)$ be the stratum of
quadratic differentials on the projective line with two simple zeros and
$6$ simple poles. Taking a double cover, ramified at each of the
poles and zeroes of the quadratic differential, yields a 
$\SL_2(\R)$-equivariant local diffeomorphism (\cite{KZ}, \cite{La})
$$\mathcal Q(1,1,-1^6) \rightarrow 
\mathcal H^\textrm{odd}(2,2)$$
The action of $\SL_2(\R)$ on the hypersurface of 
fixed area in any stratum
of quadratic differentials is ergodic by \cite{Ve2}, see also \cite{MaSm}.
These two statements imply the Lemma. 
\end{proof}
\par
Given a basis of the relative homology of a surface in $\HHH^{\rm odd}(2,2)$,
integration of the one-form defines a map to $\C^7$.
This map is a local biholomorphism (\cite{DH}, \cite{MaSm}, \cite{Ve3}), the
coordinates are called {\em period coordinates}. In the case of 
the adjusted Arnoux-Yoccoz
surface a basis of the relative homology is given 
by $\{\gamma'_2, \gamma'_3, \gamma'_6,
\gamma'_5, \beta'_1, \delta'_1, \delta'_2\}$. If $\gamma'_2=\gamma'_5$ then
the $180$ degree rotation around the center of $C'_1$ defines an
affine diffeomorphism with $8$ fixed points. Hence this equation
singles out the hyperelliptic locus $\LLL$ inside 
$\HHH^{\rm odd}(2,2)$. 
\par
We deduce from Lemma~\ref{Ratnerapplies2} that the orbit closure of
the Arnoux-Yoccoz disc contains a full neighborhood of
the initial value for $|\gamma'_3|/ \sqrt{{\rm area} (T'_1)}$, $ 
|\delta'_1|/\sqrt{{\rm area} (T'_1)}$, 
$|\gamma'_6|/\sqrt{{\rm area} (T'_2)}$ and $|\delta'_2|
/\sqrt{{\rm area} (T'_2)}$
while the other parameters are kept fixed. Using the
$\GL^+_2(\R)$-action we may vary $|\gamma'_2|/\sqrt{{\rm area} (C'_1)}$, 
$|\beta'_1|/\sqrt{{\rm area} (C'_1)}$ and ${\rm area} (C'_1)$
arbitrarily while the ratios of the areas of the splitting pieces
are kept fixed.

Let us count the dimensions that we already obtained from Lemma~\ref{Ratnerapplies2} and by applying  $\GL^+_2(\R)$. The Lie algebra $\frach$ gives 7 dimensions. We check using $\GL^+_2(\R)$ that one gets 3 more dimensions which makes 10 in total. We will now see that the ratios 
${\rm area}(T_i')/{\rm area}(C_1'),\ i = 1,2$ will give the two missing parameters.
\par
Hence till now we know that in a neighborhood of $(X',\omega')$
$$\ol{\GL^+_2(\R) \cdot (X',\omega')} \supset
\{{\rm area}(T'_i)/{\rm area}(C'_1) = \kappa_i, i=1,2\} $$
for some constants $\kappa_1$ and $\kappa_2$. We have to vary
these constants next using resplittings. The final argument
below will use yet another coordinate system. We have to
control when the hypotheses of Lemma~\ref{Ratnerapplies2} are
satisfied.
\par
Apply a  Dehn twist to the vertical cylinder of
height $(h_b)_x := (\beta_1)_x  = 1$ in Figure~\ref{fig:adjusted}. 
The corresponding direction is also a $2T2C$-direction.
We denote objects in this direction by double-primes ,
i.e.\ we have a {\em twisted splitting}
$$(X',\omega') = T''_1 \,\# \,C''_1 \,\# T''_2 \,\#\, C''_2.
$$
\par
By Lemma~\ref{Ratnerapplies2}, the $N$--orbit closure of 
$(T'_1 \#\,C'_1 \,\# T'_2 \# C'_2)$ contains
 $N\times N\times \{id\} \cdot (T'_1 \#\,C'_1 \,\# T'_2 \# C'_2)$.
In particular for $(u_1,u_2)$ in a neighborhood of $(0,0)$  
$$(X'(u_1,u_2),\omega'(u_1,u_2)) := M_{u_1}(T'_1) \,\# \,M_{u_2}C'_1 \,\#\, 
T'_2 \,\#\, M_{u_2}C'_2,$$
lies in the $N$--orbit closure, 
where $M_{u_i} = \left(\begin{array}{cc} 1 & {u_i} \\ 0 & 1 \\ \end{array}
\right)$.
\par
\begin{Lemma}
For $(u_1,u_2)$ in a neighborhood of $(0,0)$ the twisted decomposition persists.
\end{Lemma}
\par
The statement is obvious from the construction via Dehn twists.
We denote this decomposition as follows:
$$(X'(u_1,u_2),\omega'(u_1,u_2)) := T''_1(u_1,u_2) \,\# 
\,C''_1(u_1,u_2) \,\#\, 
T''_2(u_1,u_2)) \#\, C''_2(u_1,u_2).$$
We fix some more notation: $T''_i(u_1,u_2) = \R^2/\Lambda''_i(u_1,u_2)$
and $\Lambda''_3(u_1,u_2):= \Lambda''_C(u_1,u_2)$ is defined by
$C''_1(u_1,u_2) = \R^2/\Lambda''_1(u_1,u_2)$. By definition we have
$$\begin{array}{ccl}
\Lambda''_1(u_1,u_2) &= & \left\langle
\begin{pmatrix} u_1 |\gamma'_3| \\ |\gamma'_3|\end{pmatrix}, 
V(u_1,u_2) + \begin{pmatrix} (h'_s)_x + u_1 (h'_s)_y \\ (h'_s)_y
\end{pmatrix} \right\rangle\\
\Lambda''_2(u_1,u_2) &= & \left\langle
\begin{pmatrix} 0 \\ |\gamma'_6|
\end{pmatrix},  \ V(u_1,u_2) + \begin{pmatrix} (h'_m)_x  \\ (h'_m)_y
\end{pmatrix}
\right\rangle \\
\Lambda''_3(u_1,u_2) &= & \left\langle
\begin{pmatrix} u_2 |\gamma'_2| \\ |\gamma'_2| \end{pmatrix}, V(u_1,u_2) 
\right\rangle\\
\end{array}$$
where $$V(u_1,u_2) = \begin{pmatrix} 2 u_2 \ |\gamma'_2|  \\ 2|\gamma'_2|
\end{pmatrix} + \begin{pmatrix} 0  \\ |\gamma'_6|
\end{pmatrix} + \begin{pmatrix} u_1\  |\gamma'_3|  \\ |\gamma'_3|
\end{pmatrix} + \begin{pmatrix}\beta'_1  \\ 0
\end{pmatrix}.$$
\par
The next lemma shows that the $(u_1,u_2)$-twisting can indeed be used
to adjust the areas.
\par
\begin{Lemma}\label{areas}
 The map
$$\varphi: 
(u_1,u_2) \mapsto ({\rm area}(T''_1(u_1,u_2))/{\rm area}(C''_1(u_1,u_2)),
{\rm area}(T''_2(u_1,u_2))/{\rm area}(C''_1(u_1,u_2))) $$
is an invertible function in a neighborhood of $(u_1,u_2)=(0,0)$.
\end{Lemma}
\par
\begin{proof}
We remark that $\varphi$ is the composition of 
$$\psi: (u_1,u_2) \mapsto ({\rm area}(T''_1(u_1,u_2)),
{\rm area}(T''_2(u_1,u_2))) $$
and 
$$\eta: (x,y) \mapsto (\frac{2x}{1-x-y}, \frac{2y}{1-x-y}).$$
A direct computation shows that the Jacobian of $\eta$ is 
$\frac{4}{(1-x-y)^3}$. This number is non-zero when $x+y$ is far 
from 1, a condition which is satisfied for $(u_1,u_2)$ 
in a neighborhood of $(0,0)$.
\par
The rest of the proof of the lemma consists in computing the 
Jacobian of $\psi$.
\par
$${\rm area}(T''_1) = u_1 |\gamma'_3|V(u_1,u_2)_y - |\gamma'_3|
(2u_2|\gamma'_2|+u_1|\gamma'_3|+|\beta'_1|+(h'_s)_x)$$ and 
$${\rm area}(T''_2) = |\gamma'_6|(2u_2|\gamma'_2|+u_1|\gamma'_3|+
|\beta'_1|+(h'_m)_x).$$
A direct computation leads to 
$$\textrm{Jacobian}(\psi) = 2 |\gamma'_3||\gamma'_2||\gamma'_6| V_y.$$
This number is non-zero, therefore $\psi$ is locally one to one. 
\end{proof}
\par
We want to apply the result of Corollary \ref{incomm2} to the decomposition 
$(T''_1(u_1,u_2),$ $C''_1(u_1,u_2),$ 
$T''_2(u_1,u_2)) )$. The ``horizontal direction" is the direction of the vector $V(u_1,u_2)$.
We denote by $N_{V}$ the conjugate of $N$ fixing $V(u_1,u_2)$.

\begin{Lemma} \label{twisted-orbits}
 For almost all $(u_1,u_2)$ in a neighborhood of $(0,0)$ with respect
to the Lebesgue measure the $N_{V}$--orbit of 
$(T''_1(u_1,u_2), T''_2(u_1,u_2), C''_1(u_1,u_2))$
contains 
$$(G\times G \times N_V) \cdot (T''_1(u_1,u_2), T''_2(u_1,u_2), 
C''_1(u_1,u_2)).$$
\end{Lemma}
\begin{proof}

 We have to check the hypothesis of Corollary \ref{incomm2} for almost all 
 $(u_1,u_2)$ in a neighborhood of $(0,0)$ replacing the horizontal 
direction by the direction of $V(u_1,u_2)$.
 
If $\Lambda''_1(u_1,u_2)$ has a vector parallel to $V(u_1,u_2)$
there exist $t \in \R$, $n, p \in \Z$ such that
$$t V(u_1,u_2) = n \begin{pmatrix} u_1 |\gamma'_3| \\ |\gamma'_3|
\end{pmatrix} + p \begin{pmatrix} (h'_s)_x + u_1 (h'_s)_y \\ (h'_s)_y
\end{pmatrix}$$
The second coordinate equation implies 
$t = \frac{p(h'_s)_y+ n|\gamma'_3|}{2 |\gamma'_2|+ |\gamma'_6|+ 
|\gamma'_3|}.$ Thus, 
if $n$ and $p$ are fixed,  the first coordinate yields a 
linear equation in $(u_1,u_2)$.
Almost all parameters $(u_1, u_2)$, do not belong to this countable union of lines which gives the first claim. The same reasoning holds for the lattice  $\Lambda''_2(u_1,u_2)$. \\

Let $\hat{M}_t$ for $t \in \R$ denote the elements of $N_V$, 
define $A''_i := {\rm area}({\Lambda''_i}(u_1,u_2))$ 
and abbreviate $V:=V(u_1,u_2)$.
If the normalized lattice $\hat{M}_t(\widetilde{\Lambda''_2}(u_1,u_2))$ 
is commensurable to $\widetilde{\Lambda''_1}(u_1,u_2)$ for some $t$, 
there are integers $m \neq 0$, $n$, $p$ such that:
$$ 
\frac{m}{\sqrt{A_1}}\left(tV + 
\begin{pmatrix} 0 \\ |\gamma'_6|\end{pmatrix}\right)  =
\frac{1}{\sqrt{A_2}}\left(n\begin{pmatrix} u_1 |\gamma'_3| \\ |\gamma'_3|
\end{pmatrix} + p\left(V + \begin{pmatrix} (h'_s)_x + u_1 (h'_s)_y \\ (h'_s)_y
\end{pmatrix}\right)\right)
$$
We want to show that for each fixed $(m,n,p) \in \Z^3$ but arbitrary 
$t$ the set of solutions $(u_1,u_2)$ for this equation is of measure zero. 
Let $|\gamma'|:= 2|\gamma'_2| + |\gamma'_3| + |\gamma'_6|$ be the 
height of the big vertical cylinder. From the second coordinate we
obtain
$$ t =  \frac{\sqrt{A_1}(n |\gamma'_3| + p(|\gamma'| +
(h'_s)_y))}{\sqrt{A_2}m|\gamma'|} - m |\gamma'_6|.$$ 
Plugging this into the first coordinate and using the
area expressions from the proof of Lemma~\ref{areas} one obtains an algebraic
equation for $(u_1,u_2)$, whose non-triviality we have to decide.
Suppose this equation is trivial. Then $(u_1,u_2) = (0,0)$ is
a solution. But for $(u_1,u_2) = (0,0)$ we have by the Dehn
twist construction obviously
${\rm area}({\Lambda''_i}(0,0)) = {\rm area}({\Lambda'_i})$
for $i=1,2$. The field extension argument from 
Lemma~\ref{Ratnerapplies2} and Claim~\ref{extension} applies and yields
a contradiction.
\end{proof}
\par
Now, we have all the material to complete the proof of Theorem~\ref{orbit-closure}.
The surfaces $(X'(u_1,u_2))$ belong to the closure of the $\SL_2(\R)$-orbit 
of $(X', \omega')$ if $(u_1,u_2)$ belongs to a small neighborhood of $(0,0)$.
Take $(u_1,u_2)$ in  the set of full measure which satisfies  
Lemma~\ref{twisted-orbits}. The $N_V$--orbit closure of 
$(\Lambda''_1(u_1,u_2), \Lambda''_2(u_1,u_2), \Lambda''_3(u_1,u_2))$
equals 
$$
(G\times G \times N) \cdot (\Lambda''_1(u_1,u_2), 
\Lambda''_2(u_1,u_2), \Lambda''_3(u_1,u_2)).
$$ 
Applying an element of $\SL_2(\R)$, we can map $V(u_1,u_2)$ to any 
vector in $\R^2$. The previous analysis is $\SL_2(\R)$-equivariant.
Thus, given $A \in \SL_2(\R)$, the orbit closure under the action
of the unipotent subgroup of the decomposition 
of $A\cdot (X'(u_1,u_2))$ in the direction $A\cdot V$ equals the 
$G\times G \times N_{A.V}$-orbit. As in
the motivating discussion after Lemma~\ref{ergodo} we see, using
period coordinates and still for fixed $(u_1,u_2)$, that the ratios
of the splitting pieces are unchanged by this process. But 
within the real codimension $2$ neighborhood of $(X',\omega')$
in $\LLL_\alpha$ determined by the ratios
of the splitting pieces the $\SL_2(\R)$-orbit contains an open set.
On the other hand, by Lemma~\ref{areas} the data of  $(u_1,u_2)$ is 
equivalent to the areas of the splitting
pieces. Thus, the $\SL_2(\R)$-orbit closure of $(X', \omega')$ contains 
a set of positive measure in $\LLL_\alpha$. Applying 
Lemma~\ref{ergodo}, we conclude that the $\SL_2(\R)$-orbit closure 
of $(X, \omega)$ is equal to $\LLL_\alpha$. 

\appendix

\section*{Appendix}
\setcounter{section}{1}
We give in this section the details of the calculations of Theorem~\ref{pseudo-Anosov}.

We find two completely directions on the Arnoux-Yoccoz surface, namely $\theta$ and $\theta'$. 
Let us check that the pairs of slopes
$$
(\theta,\theta^\perp)=(1- \alpha^2,1+\alpha^2), 
\qquad (\theta',\theta'^\perp)=(3 + \alpha^2,\frac1{169}\left(
367+252\alpha+175\alpha^2\right))
$$
satisfy to the theorem. \medskip

Let us first introduce some notations. According to Figure~\ref{fig:gluings}
we denote the white points by $P_1,\ P_2,\ P_3$ and the black points by $P_4,\ P_5,\ P_6$ in the
order presented. The reason for 
this order will become clear in the sequel. Note that the 
points $\{P_i\}_{i=1,2,3}$ and $\{P_i\}_{i=4,5,6}$ respectively
define a singularity in $X$. If $\theta$ is a completely periodic direction 
then for each of the two singularities there are three emanating 
saddle connections. We label the saddle connection emanating from $P_i$ by $\gamma_i$.

With these conventions the labels of the singularities appear in the
counterclockwise order around a singularity. The length of $\gamma_i$ is 
denoted by $|\gamma_i| = |\int_{\gamma_i} \omega|$. We will frequently
abuse notation and write $\gamma_i$ for the vector
$
\left(\begin{smallmatrix} \int_{\gamma_i} \operatorname{Re}(\omega)  \\
\int_{\gamma_i} \operatorname{Im}(\omega) \end{smallmatrix}\right)
$
For a vector $u=(u_1,\dots,u_n)$ we will use the notation $u^2$ for $(u_1^2,\dots,u_n^2)$.

\subsection*{A first direction of slope $\theta = 1-\alpha^2$}
\label{firstdir}
\begin{figure}[htbp]
\begin{center}

\includegraphics[width=8.5cm]{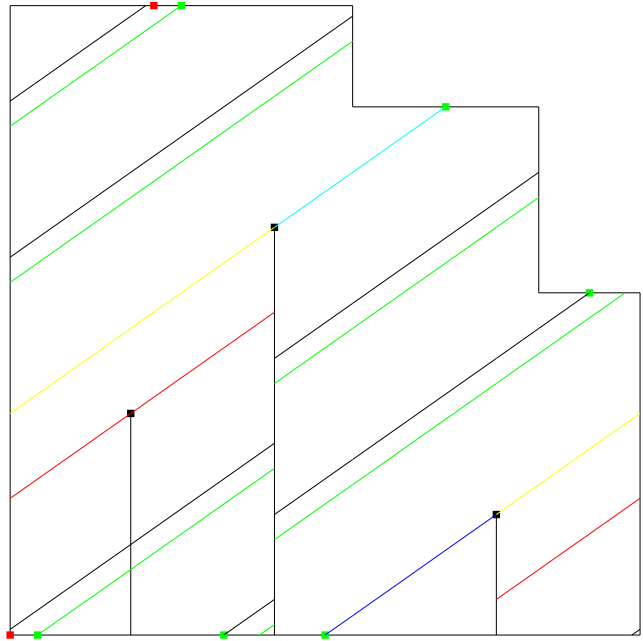}
\end{center}
\caption{\label{fig:ay:1}
Singular geodesics on the Arnoux-Yoccoz flat surface in the 
direction of slope $\theta=1-\alpha^2$: a 3C-direction.
}
\end{figure}

In this direction the surface $X$ decomposes into three metric 
cylinders. See Figure~\ref{fig:ay:1} and also Figure~\ref{fig:adjusted} (in
the vertical direction) for
the way the $6$ saddle connections bound the cylinders. We
call such a direction a {\em 3C-direction}.
A straightforward computation gives the lengths of the $\gamma_i$
(in the direction $\theta$) summarized in the following table.
$$
\begin{array}{|c|l|c|}
\hline
\textrm{saddle connection } \gamma_i & (|\gamma_i|^2)_{i=1,\dots,6}
\in \mathbb R^6 & L(\theta)^2 \in \mathbb{RP}(6) \\
\hline
i=1 & 16+18\alpha +10\alpha^2 & 40+34\alpha +22\alpha^2 \\
i=2 & 4 -6\alpha - \alpha^2 & 1 \\
i=3 & 16+ 18\alpha + 10\alpha^2 & 40+34\alpha +22\alpha^2 \\
i=4 & 4-6 \alpha + 6\alpha^2 & 4+2\alpha +2\alpha^2 \\
i=5 & 4-6\alpha - \alpha^2 & 1 \\
i=6 & 4-6\alpha + 6\alpha^2 & 4+2\alpha +2\alpha^2 \\ 
\hline
\end{array}
$$

Let us now describe the combinatorics of the 3C-direction
$\theta$. We will denote these cylinders by 
$\CCC_b, \CCC_m, \CCC_s$ for big, medium and small according
to their widths.
For the top boundaries of the cylinders, the saddle connections
appear in the following order
$$
\sigma^t_{big} = (2\ 4\ 5\ 1),\ \sigma^t_{med}=(6),\ \sigma^t_{small}=(3).
$$
And for the bottom boundaries one has
$$
\sigma^b_{big}=(3\ 2\ 6\ 5),\ \sigma^b_{med} = (4),\ \sigma^b_{small}=(1).
$$
Thus we get the following combinatorics, on which we have chosen
arbitrarily a marking.
$$
(\comb,\marking) = \left( \ 
\left(\begin{smallmatrix} \underline{2}&4&5&1 \\ \underline{3}&2&6&5\end{smallmatrix}\right),\ 
\left(\begin{smallmatrix} \underline{6} \\ \underline{4} \end{smallmatrix}\right),\ 
\left(\begin{smallmatrix} \underline{3} \\ \underline{1}\end{smallmatrix}\right)\ 
 \right)
$$
The affine holonomy vector of the big cylinder is
$\left(\begin{smallmatrix} 
4+4\alpha +2\alpha^2  \\ 
2+4\alpha +2\alpha^2\end{smallmatrix}\right).$
The 
heights of the cylinders are summarized below.
$$
\begin{array}{|c|c|c|}
\hline
\textrm{cylinder} & (\ (h_b)_y,(h_m)_y,(h_s)_y\ ) \in \mathbb R^3 &
[\ (h_b)_y,(h_m)_y,(h_s)_y\ ] \in \mathbb{RP}(3) \\
\hline
\CCC_b & 2-4\alpha +2\alpha^2 & 3+3\alpha +2\alpha^2 \\
\CCC_m & -2+2\alpha + 4\alpha^2 & 2+2\alpha +\alpha^2 \\
\CCC_s & -2+6\alpha -4\alpha^2 & 1  \\
\hline
\end{array}
$$
Equipped with this marking one can calculate the twists in the medium
and small cylinder (recall that we normalize the twist in the first,
i.e.\ the big
cylinder to zero). We first calculate the direction $\theta^\perp$.

The twist vector along the cylinder $\CCC_b$ is a vector (contained
in $\CCC_b$) joining the origin of the saddle 
connection $\gamma_3$ to the origin of the saddle connection
$\gamma_2$. This vector can be calculated in the following way. 
Let $P$ be the endpoint of the vector $h_b=\left(\begin{smallmatrix} 0 \\ 
(h_b)_y \end{smallmatrix}\right)$ based at the origin of 
$\gamma_3$, namely $P_3$. The point $P$ is located on the top of $\CCC_b$. 
Let $v_b$ be the vector joining $P$ to the origin of the
saddle connection $\gamma_2$ (namely $P_2$) in the direction
$\theta$. The twist vector of $\CCC_b$ is then $h_b + v_b$.
A simple computation shows that $v_b=\left(\begin{smallmatrix}  1-\alpha^2
\\ 2\alpha-2\alpha^2 \end{smallmatrix}\right)$. Hence the twist
vector is $h_b + v_b = \left(\begin{smallmatrix}  1-\alpha^2 \\
2-2\alpha \end{smallmatrix}\right)$. Note that here we require that 
$|v_b|<|w_b|$ (we choose the positive one with respect to the direction $\theta$). 
Finally with our marking, one obtains  
$$\theta^\perp:=\theta^\perp(0) = h_b + v_b  = \frac{2-2\alpha}{1-\alpha^2}=1+\alpha^2.$$ 
\medskip
The normalizing matrix is: 
$$M = \left(\begin{smallmatrix} 
4+4\alpha +2\alpha^2  &   & 1-\alpha^2\\ 
2+4\alpha +2\alpha^2 &  & 2-2\alpha
\end{smallmatrix}\right)^{-1} = 
\frac1{4}\left(\begin{smallmatrix} 
1+\alpha^2  &   & -4 -4\alpha-2\alpha^2\\ 
-1 &  & 6+5\alpha +3\alpha^2
\end{smallmatrix}\right).$$
Furthermore, this matrix sends the big cylinder in the direction $\theta$ to a unit square with horizontal and vertical sides.
\medskip

Let us compute the twists of cylinders $\CCC_m$ and $\CCC_s$
with respect to the direction $\theta^\perp$. We first begin 
with cylinder $\CCC_m$. 

As above let $P$ be the endpoint of the vector 
$h_m=\left(\begin{smallmatrix} 0 \\ (h_m)_y \end{smallmatrix}\right)$
based at the origin of $\gamma_4$, namely $P_4$. 
The point $P$ is located on the top of $\CCC_m$. Let
$v_m$ be the vector joining $P$ 
to the origin of the saddle connection $\gamma_6$ 
(namely $P_6$) in the direction $\theta$. 
A simple computation shows $v_m=\left(\begin{smallmatrix}
1-\alpha \\ 2-2\alpha-2\alpha^2 \end{smallmatrix}\right)$.
The vector $v_m$ is well defined up to an additive constant $n w_m$
where $n\in \mathbb Z$. 
The twist $t_m$ of $\CCC_m$ will be the difference
$$
t_m = v_m - \frac{(h_m)_y}{(h_b)_y} v_b.
$$
Thus
$$
\frac{(h_m)_y}{(h_b)_y}v_b = \frac1{2}(1+\alpha^2)
\left(\begin{smallmatrix}  1-\alpha^2 \\ 2\alpha-2\alpha^2
\end{smallmatrix}\right) = \left(\begin{smallmatrix}  1-\alpha \\
  2-2\alpha-2\alpha^2 \end{smallmatrix}\right) = v_m.
$$
Therefore the (normalized) twist $\frac{|t_m|}{|w_m|} \in [0,1[$ 
of $\CCC_m$ is zero. \medskip

Now let us finish the calculation of the twist of $\CCC_s$. 
As above let $P$ be the endpoint of the vector 
$h_s=\left(\begin{smallmatrix} 0 \\ (h_s)_y \end{smallmatrix}\right)$
based at the origin of $\gamma_1$, namely $P_1$. 
The point $P$ is located on the top of $\CCC_s$. Let
$v_s$ be the vector joining $P$ to the origin of the
saddle connection $\gamma_3$ (namely $P_3$) in the direction
$\theta$. A simple computation shows $v_s=\left(\begin{smallmatrix}
1+\alpha+2\alpha^2 \\ 2-2\alpha+2\alpha^2 \end{smallmatrix}\right)$. 
The vector $v_s$ is well defined up to an additive constant $n w_s$
where $n\in \mathbb Z$. The twist $t_s$ of $\CCC_s$ 
will be the difference $t_s = v_s -\frac{(h_s)_y}{(h_b)_y}v_b$. But
$$
\frac{(h_s)_y}{(h_b)_y}v_b = \frac1{2}(-1+2\alpha+\alpha^2)
\left(\begin{smallmatrix}  1-\alpha^2 \\ 2\alpha-2\alpha^2
\end{smallmatrix}\right) = \left(\begin{smallmatrix}  -1+\alpha+2\alpha^2 \\
-2\alpha+4\alpha^2 \end{smallmatrix}\right).
$$
Therefore the twist vector in the small cylinder is
$t_s=\left(\begin{smallmatrix} 2 \\ 2-2\alpha^2 \end{smallmatrix}\right)$.
Using $w_s =\gamma_3=\gamma_1$ we calculate the normalized twist
of $\CCC_s$ to be
$$
\frac{|t_s|}{|w_s|} =
\left(\frac{4+8\alpha-8\alpha^2}{16+18\alpha+10\alpha^2}\right)^{1/2} =
\left(1-2\alpha+\alpha^2\right)^{1/2} \ \in [0,1[
$$
We summarize the above computations for the direction $\theta=1-\alpha^2$.
$$
\begin{array}{ll}
& \\
\overrightarrow{L}^2(\theta) = 
\left[\begin{array}{c}
40+34\alpha +22\alpha^2 \\
1 \\
40+34\alpha +22\alpha^2 \\
4+2\alpha +2\alpha^2 \\
1 \\
4+2\alpha +2\alpha^2
\end{array}\right] &
\overrightarrow{W}(\theta)=
\left[\begin{array}{c}
3+3\alpha +2\alpha^2 \\
2+2\alpha +\alpha^2\\
1
\end{array} \right] \\ 
& \\
& \\
(\comb,\marking) = \left( \ 
\left(\begin{smallmatrix} \underline{2}&4&5&1 \\ \underline{3}&2&6&5\end{smallmatrix}\right),\ 
\left(\begin{smallmatrix} \underline{6} \\ \underline{4} \end{smallmatrix}\right),\ 
\left(\begin{smallmatrix} \underline{3} \\ \underline{1}\end{smallmatrix}\right)\ 
 \right) &  
\overrightarrow{T}^2(\theta,\marking,0)=
\left( 0,1-2\alpha+\alpha^2    \right) \\
& \textrm{normalized twists} \\
& \\
\end{array}
$$

\subsection*{A second direction of slope $\theta' = 3+\alpha^2$} \label{seconddir}

\begin{figure}
\begin{center}

\includegraphics[width=8.5cm]{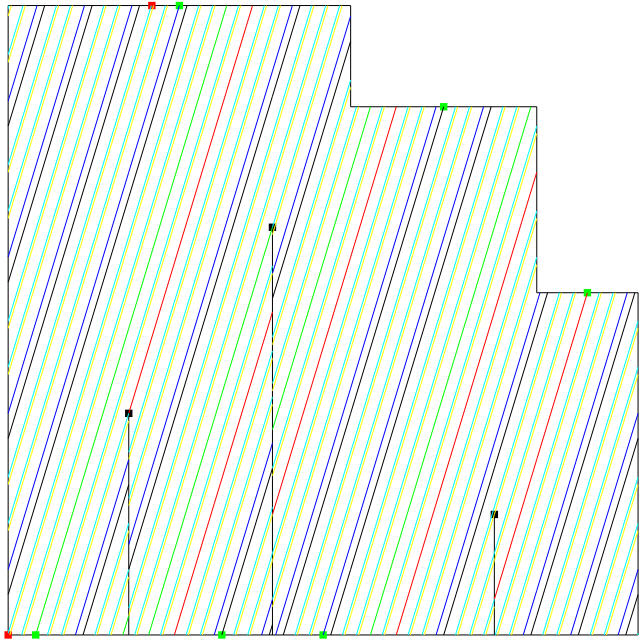}
\end{center}
\caption{\label{fig:ay:2}
Singular geodesics on the Arnoux-Yoccoz flat surface in the 
direction $\theta'=3+\alpha^2$: again a 3C-direction.
}
\end{figure}

Using Figure~\ref{fig:ay:2} one checks that $\theta'$ is also
a 3C-direction with the following invariants.
$$
\begin{array}{|c|c|c|}
\hline
\textrm{saddle connection } \gamma'_i & (|\gamma'_i|^2)_{i=1,\dots,6}
\in \mathbb R^6 & L(\theta')^2 \in \mathbb{RP}(6) \\
\hline
i=1 & 26+20\alpha +11\alpha^2 & 1 \\
i=2 & 144 +114\alpha + 74\alpha^2 & 4+2\alpha +2\alpha^2 \\
i=3 & 144 +114\alpha + 74\alpha^2 & 4+2\alpha +2\alpha^2 \\
i=4 & 26+20 \alpha +11\alpha^2 & 1 \\
i=5 & 1612+1354\alpha + 878\alpha^2 & 40+34\alpha +22\alpha^2 \\
i=6 & 1612+1354\alpha + 878\alpha^2 & 40+34\alpha +22\alpha^2 \\ 
\hline
\end{array}
$$
As previously we label the cylinders 
by $\CCC'_b, \CCC'_m, \CCC'_s$ according to their widths. The 
combinatorics with some choice of marking is given by 
$$
(\comb',\marking') = \left( \ 
\left(\begin{smallmatrix} \underline{4}&3&1&6 \\ \underline{5}&4&2&1\end{smallmatrix}\right),\ 
\left(\begin{smallmatrix} \underline{2} \\ \underline{3} \end{smallmatrix}\right),\ 
\left(\begin{smallmatrix} \underline{5} \\ \underline{6}\end{smallmatrix}\right)\ 
 \right).
$$
A straightforward computation shows the expected results
for the heights of the cylinders in the direction $\theta'$.
$$
\begin{array}{|c|c|c|}
\hline
\textrm{cylinder} & (\ (h'_b)_y,(h'_m)_y,(h'_s)_y\ ) \in \mathbb R^3 &
[\ (h'_b)_y,(h'_m)_y,(h'_s)_y \ ] \in \mathbb{RP}(3) \\
\hline
\CCC'_b & -6+8\alpha +6\alpha^2 & 3+3\alpha +2\alpha^2 \\
\CCC'_m & -2+6\alpha - 4\alpha^2 & 2+2\alpha +\alpha^2 \\
\CCC'_s & 10-14\alpha -8\alpha^2 & 1 \\
\hline
\end{array}
$$
To calculate the twists of $\CCC'_m$ 
and $\CCC'_s$ we first calculate the direction $\theta'^\perp$.
The twist vector of $\CCC'_b$ is a vector (contained in this
cylinder) joining the origin of the saddle $\gamma'_5$ to the origin of
the saddle connection $\gamma'_4$. We will use this remark for the
computation of the twist vector of $\CCC'_b$.

Let $P$ be the endpoint of the vector $h'_b=\left(\begin{smallmatrix} 0
\\ (h'_b)_y \end{smallmatrix}\right)$ based at the origin of
$\gamma'_5$, namely $P_5$. The point $P$ is located on the top of $\CCC'_b$. 
Let $v'_b$ be the vector joining $P$ to the origin of the
saddle connection $\gamma'_4$ (namely $P_4$) in the direction
$\theta'$. A simple computation  using the
requirement $|v'_b|<|w'_b|$ shows $v'_b=\left(\begin{smallmatrix}
5+2\alpha-\alpha^2 \\ 18+2\alpha \end{smallmatrix}\right)$.
Given $n'_0 \in \mathbb Z$, the twist vector of $\CCC'_b$ is 
then $h'_b + v'_b + n'_0 w'_b$. 
The width of the big cylinder is given by the vector
$w'_b = \gamma'_5 + \gamma'_4 + \gamma'_2 + \gamma'_1 =  
\left(\begin{smallmatrix}  14+12\alpha+8\alpha^2 \\
46+40\alpha+26\alpha^2 \end{smallmatrix}\right)$. 
Finally with our marking, one obtains the direction 
$\theta'^\perp := \theta'^\perp(n_0) =h'_b + v'_b + n'_0 w'_b$.  
\par
\medskip \noindent
Proceeding as in the previous section we obtain
$v'_m=\left(\begin{smallmatrix}
4-3\alpha-\alpha^2 \\ 10-8\alpha+4\alpha^2 \end{smallmatrix}\right)$
and 
\begin{multline*}
\frac{(h'_m)_y}{(h'_b)_y}(v'_b + n'_0 w'_b) = \frac1{2}(1+\alpha^2)
\left(
\left(\begin{smallmatrix}  5+2\alpha-\alpha^2 \\ 18+2\alpha
\end{smallmatrix}\right)+ 
n'_0 \left(\begin{smallmatrix}  14+12\alpha+8\alpha^2 \\ 46+40\alpha+26\alpha^2
\end{smallmatrix}\right)
\right)  = \\
= \left(\begin{smallmatrix}  4-\alpha+\alpha^2 \\
  10+8\alpha^2 \end{smallmatrix}\right) + 
n'_0 \left(\begin{smallmatrix}  9+8\alpha+5\alpha^2 \\
  30+26\alpha+16\alpha^2 \end{smallmatrix}\right).
\end{multline*}
This yields
$$t'_m = v'_m -\frac{(h'_m)_y}{(h'_b)_y}(v'_b+n'_0 w'_b)
= \left(\begin{smallmatrix}  -2\alpha-2\alpha^2 \\
  -8\alpha-4\alpha^2 \end{smallmatrix}\right) -
n'_0 \left(\begin{smallmatrix}  9+8\alpha+5\alpha^2 \\
  30+26\alpha+16\alpha^2 \end{smallmatrix}\right).
$$
The width of $\CCC'_m$ is the norm of the vector
$w'_m=\gamma_2=\gamma_3= \left(\begin{smallmatrix}  3+2\alpha+\alpha^2
\\ 10+6\alpha+4\alpha^2 \end{smallmatrix}\right)$. Therefore the normalized
twist of $\CCC'_m$ is
\begin{multline}
\label{formulae:tm}
\frac{|t'_m|}{|w'_m|} = \left(4+2n'_0+7{n'_0}^2 + (-6+6{n'_0}^2)\alpha +
(-2+2n'_0+4{n'_0}^2)\alpha^2\right)^{1/2} \ \in \mathbb R / \mathbb Z 
\end{multline}

\begin{Remark}
If $n'_0=-1$ then $\frac{|t'_m|}{|w'_m|}$ equals $3$. 
Thus in case $n'_0=-1$, the normalized twist of $\CCC'_m$ is
zero, more precisely $t'_m = 3w'_m$.
\end{Remark}
\par
Now let us compute the twist of the cylinder $\CCC'_s$. We obtain $v'_s=\left(\begin{smallmatrix}
-1+5\alpha+6\alpha^2 \\ -4+22\alpha+12\alpha^2
\end{smallmatrix}\right)$ and 
\begin{multline*}
\frac{(h'_s)_y}{(h'_b)_y}(v'_b+n'_0 w'_b) = 
\frac1{2}(-1+2\alpha+\alpha^2) \left(
\left(\begin{smallmatrix}  5+2\alpha-\alpha^2 \\ 18+2\alpha
\end{smallmatrix}\right) + n'_0 
\left(\begin{smallmatrix}  14+12\alpha+8\alpha^2 \\
  46+40\alpha+26\alpha^2 \end{smallmatrix}\right)
\right) =  \\
= \left(\begin{smallmatrix}  -2+3\alpha+5\alpha^2 \\
  -8+16\alpha+10\alpha^2 \end{smallmatrix}\right) + n'_0
\left(\begin{smallmatrix}  3+2\alpha+\alpha^2 \\
  10+6\alpha+4\alpha^2 \end{smallmatrix}\right)
\end{multline*}
Together this gives
$$t'_s = v'_s -\frac{(h'_s)_y}{(h'_b)_y}(v'_b+ n'_0 w'_b)
= \left(\begin{smallmatrix} 1+2\alpha+\alpha^2 \\
4+6\alpha+2\alpha^2 \end{smallmatrix}\right) -n'_0 
\left(\begin{smallmatrix}  3+2\alpha+\alpha^2 \\
10+6\alpha+4\alpha^2 \end{smallmatrix}\right).
$$
The width $w'_s$ of the small cylinder is given by 
$w_s =\gamma'_5=\gamma'_6$, hence
\begin{multline}
\label{formulae:ts}
\frac{|t'_s|}{|w'_s|} = \left(2-{n'_0}^2 + (-2+2n'_0+{2n'_0}^2)\alpha +
(-3-4n'_0)\alpha^2\right)^{1/2} \ \in \mathbb R / \mathbb Z
\end{multline}

\begin{Remark}
If $n'_0=-1$ then $\frac{|t'_s|}{|w'_s|}$ equals 
$(1-2\alpha+\alpha^2)^{1/2} \in [0,1[$.
\end{Remark}
\par
One can summarize the above computations (for the direction 
$\theta'=3+\alpha^2$)
$$
\begin{array}{ll}
& \\
\overrightarrow{L}^2(\theta') = 
\left[\begin{array}{c}
1 \\
4+2\alpha +2\alpha^2 \\
4+2\alpha +2\alpha^2 \\
1 \\
40+34\alpha +22\alpha^2 \\
40+34\alpha +22\alpha^2 \\
\end{array}\right] &
\overrightarrow{W}(\theta')=
\left[\begin{array}{c}
3+3\alpha +2\alpha^2 \\
2+2\alpha +\alpha^2\\
1
\end{array} \right] \\ 
& \\
& \\
(\comb',\marking') = \left( \ 
\left(\begin{smallmatrix} \underline{4}&3&1&6 \\ \underline{5}&4&2&1\end{smallmatrix}\right),\ 
\left(\begin{smallmatrix} \underline{2} \\ \underline{3} \end{smallmatrix}\right),\ 
\left(\begin{smallmatrix} \underline{5} \\ \underline{6}\end{smallmatrix}\right)\ 
 \right) & 
\overrightarrow{T}^2(\theta',\marking',n'_0)= 
(\ \textrm{see}~(\ref{formulae:tm}),\ \textrm{see}~(\ref{formulae:ts})\ ) \\
& \textrm{normalized twists}  \\
& \\
\end{array}
$$
The normalizing matrix is
$$M' = \left(\begin{smallmatrix} 
14+12\alpha +8\alpha^2  &   & -9-10\alpha-9\alpha^2\\ 
46+40\alpha +26\alpha^2 &  & -34-30\alpha-20\alpha^2
\end{smallmatrix}\right)^{-1} = 
\frac1{4}\left(\begin{smallmatrix} 
-49-42\alpha-27\alpha^2  &   & -66 -56\alpha-36\alpha^2\\ 
14+14\alpha+9\alpha^2 &  & 20+17\alpha +11\alpha^2
\end{smallmatrix}\right).$$
\medskip

\subsection*{A second pseudo-Anosov diffeomorphism}
\label{sec:derivative}

We have now all necessary tools to prove that there exists another 
pseudo-Anosov diffeomorphism, say $\tilde \Phi$, in 
the stabilizer of the Arnoux-Yoccoz Teichm\"uller disc.

\begin{proof}
If we take $n'_0=-1$ then 
Theorem~\ref{theo:invariants} applies with the above $\theta,\ \theta'$, 
and permutations 
$\pi_1$ on $k=6$ elements and $\pi_2$ on $p=3$ elements given by
$$
\pi_1 = \left( \begin{tabular}{cccccc}
1 & 2 & 3 & 4 & 5 & 6 \\
6 & 4 & 5 & 3 & 1 & 2
\end{tabular} \right) \qquad \textrm{and} \qquad
\pi_2 = \textrm{id}
$$
We claim that the diffeomorphism $\tilde\Phi$ thus obtained
has the required property. For this purpose we calculate
its derivative in the coordinates $(x,y)$ given by the Arnoux-Yoccoz
construction at the beginning of Section~\ref{sec:calculus}.
By construction, one knows the image by
$\tilde \Phi$ of the saddle connections in direction $\theta$.
In particular, $\tilde \Phi (\gamma_i) = \gamma'_{\pi_1(i)}$.
For  $i=3$ we obtain
$$
\tilde \Phi (\gamma_3) = \gamma'_5 \qquad \textrm{ i.e. } \qquad 
D \tilde \Phi \left(\begin{smallmatrix} 3+2\alpha+\alpha^2 \\
2+2\alpha \end{smallmatrix}\right) =
\left(\begin{smallmatrix} 9+8\alpha+5\alpha^2 \\
30+26\alpha+16\alpha^2 \end{smallmatrix}\right) 
$$
Moreover, with our previous normalization, the image of the vector $h_b + v_b$
is the vector $h'_b + v'_b + n'_0 w'_b$ (with $n'_0=-1$). Therefore
$$
D \tilde \Phi \left(\begin{smallmatrix} 1-\alpha^2 \\ 2-2\alpha
\end{smallmatrix}\right) =
\left(\begin{smallmatrix} -9-10\alpha-9\alpha^2 \\
  -34-30\alpha-20\alpha^2 \end{smallmatrix}\right).
$$
Applying a change of basis, we get the matrix $D \tilde \Phi$ in 
the coordinate $(x,y)$:
$$
D \tilde \Phi = \left(
\begin{smallmatrix} 9+8\alpha+5\alpha^2 & &-9-10\alpha-9\alpha^2 \\ \\
30+26\alpha+16\alpha^2 & & -34-30\alpha-20\alpha^2 
\end{smallmatrix}
\right)
\left(
\begin{smallmatrix} 3+2\alpha+\alpha^2 & &1-\alpha^2 \\ \\
2+2\alpha & & 2-2\alpha
\end{smallmatrix}
\right)^{-1} = \left(
\begin{smallmatrix} 23+18\alpha+12\alpha^2 & &-29-24\alpha-16\alpha^2 \\ \\
74+62\alpha+40\alpha^2 & & -95-80\alpha-52\alpha^2 
\end{smallmatrix}
\right).
$$
Obviously, $D \tilde \Phi = M'^{-1}M$. 
Since the trace of this matrix is greater than $2$, the diffeomorphism $\tilde
\Phi$ is a pseudo-Anosov diffeomorphism. 
\end{proof}


\end{document}